# On the admissibility of Horvitz-Thompson estimator for estimating causal effects under network interference*

**Vishesh Karwa and Edoardo M. Airoldi**

*Temple University*
*Department of Statistics, Operations and Data Science*

**Abstract:**
The Horvitz–Thompson (H-T) estimator is widely used for estimating network causal effects. We study its optimality properties by embedding it in the class of all linear estimators. We show that, under any form of interference, the H-T estimator is inadmissible with respect to mean squared error for designs that generate a random number of units for a given treatment-exposure combination, which includes completely randomized and Bernoulli designs. In contrast, the Hájek and Ratio estimators are admissible under such designs as long as their weights exceed one. We show that the H-T estimator becomes admissible under restricted randomization schemes that fix treatment-exposure counts.

These findings motivate a new estimator—the conditional Horvitz–Thompson (CHT) estimator—which uses propensity scores conditional on the realized number of units in each treatment–exposure combination. We prove that the CHT estimator is unbiased and admissible for network causal effects, for any design, within the class of linear estimators. Because computing conditional propensity scores can be challenging in practice, we develop three practical approximations: an MCMC approximation and two closed-form approximations based on the difference-in-means and Ratio estimators, the first one avoids computing any propensity scores. Numerical experiments on the Add Health network illustrate the theoretical results.

**Keywords and phrases:** causal inference; potential outcomes; average treatment effect; total treatment effect; network interference; Horvitz-Thompson; Hajek estimator; Ratio estimator; Admissibility.

---

*Vishesh Karwa is an Assistant Professor in the Department of Statistics, Operations and Data Science, Fox Business School at Temple University. Edoardo M. Airoldi is the Millard E. Gladfelter Professor of Statistics & Data Science and the Director of the Data Science Center at the Fox Business School, Temple University. This work was partially supported by National Science Foundation awards CAREER IIS-1149662 and IIS-1409177, and by Office of Naval Research YIP award N00014-14-1-0485, all to Harvard University.

## 1. Introduction, Summary and Related work

### *1.1. Motivation*

Estimating causal effects under treatment interference in finite populations poses three main challenges: (a) Unlike the classical setting (Neyman, 1923; Rubin, 1974; Cox, 1958), the number of potential outcomes per unit explodes because of interference and depends on the assumed interference model; (b) there are multiple, non-equivalent definitions of causal effects; and (c) the difference-in-means estimator is generally biased, see Karwa and Airoldi (2018) for more details of these issues.

Each issue has motivated distinct solutions. For example, to address (a), researchers specify an interference (or exposure) model to reduce the number of potential outcomes (Hudgens and Halloran, 2008b; Tchetgen and VanderWeele, 2012; Manski, 2013; Bowers et al., 2013; Aronow and Samii, 2017; Sussman and Airoldi, 2017; Viviano, 2020; Leung, 2022). For (b), various definitions of network causal effects (estimands) and unbiased estimators have been proposed (Sävje et al., 2021; Hu et al., 2022; Yu et al., 2022; Choi, 2023). The two main estimand classes are the *average treatment effects* and *expected average treatment effects*. Average treatment effects contrast potential outcomes at two fixed treatment-exposure combinations, whereas expected average treatment effects contrast marginal or conditional expectations of potential outcomes, treating the exposure level as a random variable under the design, see Karwa and Airoldi (2018) for a discussion of the difference between these two classes. Finally, to address (c), the Horvitz–Thompson (H-T) estimator (Horvitz and Thompson, 1952) is widely used for unbiased estimation of average treatment effects under various interference models. In this paper, we focus on estimation of average treatment effects under treatment interference and study the optimality of linear weighted estimators, including the H-T estimator.

### *1.2. Summary of contributions*

We systematically study the H-T and other linear weighted estimators for estimating *average treatment effects* under network interference across common designs. Our contributions are as follows:

1. *Unbiasedness:* The H-T estimator is known to be unbiased for average treatment effects (Aronow and Samii, 2017). We illustrate that this property depends not only on the design but also on the interference model and estimand type. While design and estimand are chosen by the experimenter, the interference form must be modeled and may be misspecified. Thus, unbiasedness of the H-T estimator requires correct model specification. We derive analytical expressions for H-T weights under various exposure models, useful for practitioners who otherwise rely on Monte Carlo simulation.
2. *Optimality:* Assuming the interference model is correct, we analyze the optimality of the H-T estimator within a broader class of weighted linear estimators. We show it is inadmissible w.r.t. mean square error for many designs. The reason lies in a feature unique to network interference: the number of units at a given exposure level is indirectly assigned, and hence random. We call such designs *random exposure designs*. This property leads to a proof of inadmissibility of the H-T estimator for estimating network average causal effects under completely randomized and Bernoulli designs. In contrast, when the number of units per exposure level is fixed, we call such designs *fixed exposure designs*, the H-T estimator is admissible.
3. *Admissible estimators:* After establishing inadmissibility of the H-T estimator for random exposure designs, we characterize a class of admissible linear estimators. We show that a linear estimator is inadmissible if its minimum variance over all potential outcomes is strictly greater than zero, and admissible if it achieves zero variance for some potential outcome vector, as long as its weights are greater than 1. Using these results, we prove that the Hajek and ratio estimators are admissible under any design as long as their weights exceed one.
4. *Conditional H-T estimator:* We propose a new estimator, the *conditional Horvitz–Thompson* estimator and show it is unbiased and admissible under any design. It uses propensity scores conditional on the number of units with a given treatment-exposure status which is an ancillary statistic. It reduces to the standard H-T estimator for fixed exposure designs. Exact computation of the conditional propensity scores can be difficult under interference, so we develop three approximations: two analytical and one empirical using MCMC. We evaluate these approximations and other results via simulations.

### ***1.3. Related Work***

The concept of exposure mappings under network interference was formalized by Manski (2013) and Aronow and Samii (2017). Tchetgen and VanderWeele (2012) introduced the H-T estimator for partial interference (Sobel, 2006), which was later extended to network interference by Aronow and Samii (2017). For expected average treatment effects, several studies show estimation is possible even under unknown interference. Sävje et al. (2021) define expected average treatment effects (EATEs) as contrasts between conditional expectations of potential outcomes (expectation w.r.t. the design) and establish conditions on interference growth for unbiased estimation by difference-of-means. Sävje (2023) show that even with misspecified exposure mappings, EATEs can be estimated by redefining the estimand as a contrast between the expectations of the potential outcomes conditional on the "mis-specified" exposure mappings. Leung (2022) relax neighborhood interference to allow $k$-hop dependence with decay, and Hu et al. (2022) propose a complementary estimand of expected average interference effect.

In contrast, estimating *average treatment effects* typically requires knowledge of the exposure model. However, Yu et al. (2022) show that the total treatment effect—an average treatment effect—can be estimated without this knowledge under a heterogeneous additive model of potential outcomes. Using individual baseline effects, they propose an estimator of the total treatment effect not requiring the exposure model. Gao and Ding (2023) provide regression-based point and standard error estimators satisfying design-based properties.

Another line of work uses hypothesis testing framework. Rosenbaum (2007) invert randomization tests to form confidence intervals under interference; Athey et al. (2017) develop exact tests for non-sharp nulls; Basse et al. (2019) propose graph-theoretic randomization tests; and Aronow (2012); Pouget-Abadie et al. (2019) test for the presence of interference. A complementary line proposes new designs for estimating *average treatment effects* (Toulis and Kao, 2013; Ugander et al., 2013; Eckles et al., 2016; Jagadeesan et al., 2020; Ugander and Yin, 2023). Hudgens and Halloran (2008a) introduce two-stage randomization for partial interference, see also Basse and Feller (2018). Toulis and Kao (2013) propose sequential randomization for peer effects, while Ugander et al. (2013) present graph cluster randomization for total treatment effects. Eckles et al. (2016) develop design and analysis methods to reduce bias due interference, and Jagadeesan et al. (2020) propose design principles ensuring the difference-in-means estimator is unbiased under additive models of interference. Basse and Airoldi (2018) propose model-assisted designs for network-correlated outcomes. Work on optimality of estimators under network interference is limited. An exception is Sussman and Airoldi (2017), who derive a minimum integrated variance estimator among all linear unbiased estimators for average treatment effects.

## 2. Potential Outcomes and Average treatment effects under Network Interference

Consider $n$ units indexed by $\{1, 2, \ldots, n\}$ and a binary treatment vector $\mathbf{z} \in \{0, 1\}^n$. Each unit has potential outcomes $\{Y_i(\mathbf{z})\}$, depending on the entire treatment vector $\mathbf{z}$. In the SUTVA setting, each unit's outcome depends only on its own treatment, i.e., $Y_i(\mathbf{z}) = Y_i(z_i)$ (Cox, 1958; Rubin, 1978). Under network interference, $Y_i(\mathbf{z})$ also depends on other units' treatments. We formalize this through the *network interference assumption*. Let $z_{N_i}$ denote the treatment vector of unit $i$'s interference neighborhood $N_i$, then $Y_i(\mathbf{z}) = Y_i(z_i, z_{N_i})$. This reduces the number of potential outcomes per unit from $2^n$ to $2^{|N_i|+1}$.

We further use the *exposure mapping assumption* (Aronow and Samii, 2017; Manski, 2013), which states that $Y_i(z_i, z_{N_i}) = Y_i(z_i, f(z_{N_i})) = Y_i(z_i, e_i)$, where $e_i = f(z_{N_i}) \in \{0, 1, \ldots, K_i - 1\}$ is the *exposure* level of unit $i$. Thus, each unit has $2K_i$ potential outcomes. Let $d_j = (z, e)$ denote a treatment–exposure combination, with $z \in \{0, 1\}$ and $e \in \{0, \ldots, K_i - 1\}$. For simplicity, assume $K_i = K$ for all $i$. We define two special exposures: $e = 0$ (no exposure) and $e = 1$ (full exposure). Note that the interference neighborhoods $N_i$ may overlap. We represent them by an *interference graph* $G$, where $N_i$ is the adjacency list of node $i$. Thus, $G$ encodes all interference relations. We assume $G$ is fixed; some works relax this via random graph models (Li and Wager, 2022; Li et al., 2021). We also make the *Consistency assumption,* as in classical causal inference, that there are no hidden versions of potential outcomes. Hence, the observed outcome is given by $Y_i^{obs} = \sum_{z,e} Y_i(z, e) I(Z_i = z, E_i = e)$.

Given the network interference and exposure mapping assumptions, we define *average treatment effects* as the average contrast between two fixed potential outcomes:

$$\theta(d_1, d_0) = \frac{1}{n} \sum_{i=1}^{n} [Y_i(d_1) - Y_i(d_0)], \tag{1}$$

where $d_1 = (z, e)$ and $d_0 = (z', e')$ are any two fixed treatment–exposure pairs. We sometimes omit the dependence on $(d_1, d_0)$ for brevity. For illustration, consider the $2 \times 2$ *exposure model* with $e_i = I(\sum_{j \in N_i} z_j > 1)$. Here $e_i = 1$ if a unit has at least one treated neighbor, and 0 otherwise. Then a non-parametric representation of the potential outcomes is

$$Y_i(z_i, e_i) = \alpha_i + \beta_i z_i + \gamma_i e_i + \delta_i z_i e_i,$$

where $z_i, e_i \in \{0, 1\}$. If $(z, e) = (1, 1)$ and $(z', e') = (0, 0)$, $\theta$ is the *total treatment effect*; if $(z, e) = (1, 0)$ and $(z', e') = (0, 0)$, the *direct treatment effect*; if $(z, e) = (0, 1)$ and $(z', e') = (0, 0)$, the *additive interference effect*; and if $(z, e) = (1, 1)$ and $(z', e') = (1, 0)$, the *total interference effect* (Karwa and Airoldi, 2018).

A design randomly assigns $\mathbf{z}$ according to a distribution $P(\mathbf{Z} = \mathbf{z})$. Two common designs are the Bernoulli and completely randomized design (CRD). In the *Bernoulli design,* each $z_i$ is independently set to 1 with probability $p$ and 0 otherwise: $P(\mathbf{Z} = \mathbf{z}) = p^{|\mathbf{z}|}(1 - p)^{n-|\mathbf{z}|}$. In the *CRD*, $n_t$ of $n$ units are assigned to treatment and the rest to control: $P(\mathbf{Z} = \mathbf{z}) = 1/\binom{n}{n_t}$, if $|\mathbf{z}| = n_t$, 0, otherwise. In both designs, the number of units assigned to an exposure level $e_i = f(z_{N_i})$ are random and depend on $N_i$. Under CRD, the number of treated units is fixed but the number in each exposure level is random, so treatment–exposure counts are random. Under the Bernoulli design, both treatment and exposure counts are random.

For any treatment–exposure pair $d = (z, e)$, define $\Omega_i(z, e) = \{\mathbf{z} : z_i = z, e_i = e\}$ as the set of assignment vectors revealing $Y_i(z, e)$. Let $D_i = (Z_i, E_i)$ be unit $i$'s realized treatment–exposure pair under design $P(\mathbf{Z})$. The corresponding propensity scores are

$$\pi_i(d_1) = P(D_i = d_1) = \sum_{\mathbf{z}\in\Omega_i(z,e)} P(\mathbf{Z} = \mathbf{z}), \qquad \pi_i(d_0) = P(D_i = d_0) \sum_{\mathbf{z}\in\Omega_i(z',e')} P(\mathbf{Z} = \mathbf{z}).$$

Each unit's propensity score depends on the design $P(\mathbf{Z})$, exposure mapping $f(z_{N_i})$, and chosen treatment–exposure pair $(z, e)$. The Horvitz–Thompson (H–T) estimator treats observed outcomes as a random sample drawn with unequal propensity scores:

$$\hat{\theta}_{HT} = \frac{1}{n}\left(\sum_{i=1}^{n} \frac{I(D_i = d_1)}{\pi_i(d_1)} Y_i(d_1) - \sum_{i=1}^{n} \frac{I(D_i = d_0)}{\pi_i(d_0)} Y_i(d_0)\right).$$

## 3. Linear estimators for Causal Effects

It is well known that the H–T estimator is unbiased for estimating average treatment effects. Here we embed it in a larger class of linear weighted estimators and characterize when unbiasedness holds. We first consider the class of linear weighted unbiased estimators, then restrict to a smaller class where the H–T estimator turns out to be the unique unbiased member.

*Linear weighted unbiased estimators:* For any design $P(\mathbf{Z} = \mathbf{z})$, unbiased estimators can be constructed using techniques from survey sampling. Following

Godambe (1955), consider the general class of linear weighted estimators

$$\hat{\theta}_1(\mathbf{z}) = \sum_i w_i(\mathbf{z}) Y_i^{obs}, \tag{2}$$

where $w_i(\mathbf{z})$ is the weight for unit $i$ under assignment $\mathbf{z}$. The weights in general depend on the entire assignment vector. The set of weights yielding unbiased estimators is characterized by a system of linear equations involving the design, interference neighborhoods, and the exposure model.

**Theorem 3.1.** *Consider an exposure model $e_i = f(\mathbf{z}_{N_i})$ where $N_i$ is the interference neighborhood of unit i. Let $\Omega_i(z, e) = \{\mathbf{z} : z_i = z, e_i = e\}$. Similarly, let $\Omega_i(z', e') = \{\mathbf{z} : z_i = z', e_i = e'\}$. The estimator $\hat{\theta}$ in equation 2 is unbiased for $\theta = \frac{1}{n}\sum_i (Y_i(z, e) - Y_i(z', e'))$ if and only if $0 < \pi_i(z, e) < 1$ and $0 < \pi_i(z', e') < 1$ and $w_i(\mathbf{z})$ satisfy the following system of equations:*

$$\begin{aligned}
&\sum_{\mathbf{z}\in\Omega_i(z,e)} w_i(\mathbf{z})P(\mathbf{Z}) = \frac{1}{n}, \quad \forall i = 1, \ldots, n, \\
&\sum_{\mathbf{z}\in\Omega_i(z',e')} w_i(\mathbf{z})P(\mathbf{Z}) = -\frac{1}{n}, \quad \forall i = 1, \ldots, n, \\
&\sum_{\mathbf{z}\in\Omega_i(z^*,e^*)} w_i(\mathbf{z})P(\mathbf{Z}) = 0, \quad \forall (z^*, e^*) \neq (z, e), (z', e'),\ i = 1, \ldots, n.
\end{aligned}$$

Theorem 3.1 (proof in the appendix) implies that, even with a fixed graph, exposure model, and design, the linear system typically has infinitely many solutions; hence there may be infinitely many unbiased estimators. Choosing the best among them by minimizing expected variance under a prior was studied by Sussman and Airoldi (2017). If we restrict the weights so that each depends on $\mathbf{z}$ only via $(z_i, e_i)$, we get the smaller class

$$\hat{\theta}_2 = \sum_i w_i(z_i, e_i) Y_i^{obs}. \tag{3}$$

This reduction that weights depend only on a unit's realized treatment–exposure pair is a form of sufficiency, and is natural since potential outcomes are also reduced from $Y_i(\mathbf{z})$ to $Y_i(z_i, e_i)$. Under this restriction, Theorem 3.2 shows the H–T estimator is the unique unbiased estimator in the class (3).

**Theorem 3.2.** *Consider the estimators of type $\hat{\theta}_2$ given by equation 3. Without any further assumptions on the potential outcomes, the only unbiased estimator of $\theta$ in* *this class is the Horvitz-Thompson estimator* $\hat{\theta}_{HT}$ *where*

$$w_i(z_i, e_i) = \begin{cases} \frac{1}{n\pi_i(z,e)} & \textit{if } (z_i, e_i) = (z, e), \\ -\frac{1}{n\pi_i(z',e')} & \textit{if } (z_i, e_i) = (z', e'), \\ 0 & \textit{otherwise}, \end{cases}$$

*where* $\pi_i(z, e) = \sum_{\mathbf{Z} \in \Omega_i(z,e)} P(\mathbf{Z})$ *and* $\pi_i(z', e') = \sum_{\mathbf{Z} \in \Omega_i(z',e')} P(\mathbf{Z})$.

Hence, if we insist on unbiasedness and weights of the reduced form $w_i(z_i, e_i)$, the H–T estimator is the only choice. It is then natural to ask whether the H–T estimator is optimal within the larger class (2); the next section studies its admissibility in the larger class of estimators. Theorem 3.2 also shows H–T weights depend on the unit, the design, and the exposure model. In the Appendix, we derive closed-form weights for common exposure models and designs which can be useful for practitioners.

## 4. Admissibility of the Horvitz-Thompson estimator for network causal effects

In this section we study the admissibility of the H–T estimator w.r.t. mean squared error (MSE) within the class of all linear estimators defined in (2). For an estimator $\hat{\theta}$, write $\text{MSE}(\hat{\theta}) = \mathbb{E}_{P(\mathbf{Z})}[(\hat{\theta} - \theta)^2]$, where the expectation is over the design.

**Definition 1** (Admissibility)**.** *For a fixed design, a linear estimator* $\hat{\theta}_1$ *is admissible if no other linear estimator* $\hat{\theta}_2$ *satisfies* $\text{MSE}(\hat{\theta}_2) \leq \text{MSE}(\hat{\theta}_1)$ *for all* $\theta$ *with strict inequality for at least one* $\theta$.

In survey sampling with fixed sample-size designs, it is well known that the H–T estimator is admissible for finite-population totals (Godambe and Joshi, 1965). Our results show that this classic result can be extended to the causal inference setting without interference. In particular, for a CRD, the H-T estimator for estimating the average causal effect without interference is admissible with respect to the MSE, in the class of all linear estimator. In contrast, for network causal effects under common designs like CRD or Bernoulli, we show that the H–T estimator is *inadmissible* among all linear estimators. The reason is that, under interference, the number of units assigned to a given treatment–exposure pair is generally random; this extra randomness can be exploited to construct estimators with uniformly smaller MSE. We call designs where the number of units assigned to a given treatment-exposure combination as *random exposure designs*.

**Definition 2** (Random exposure designs)**.** *Let $\theta$ be the contrast in* (1) *between $d_1 = (z, e)$ and $d_0 = (z', e')$. Define $n_{d_1} = \sum_{i=1}^{n} I(Z_i = z, E_i = e)$ and $n_{d_0} = \sum_{i=1}^{n} I(Z_i = z', E_i = e')$. A design $P(\mathbf{Z})$ is a* random exposure design *for $\theta$ if* $\mathrm{Var}(n_{d_1}) > 0$ *and* $\mathrm{Var}(n_{d_0}) > 0$.

**Theorem 4.1** (Inadmissibility of H–T)**.** *Let $P(\mathbf{Z})$ be any random exposure design as in Definition 2. Let $\theta$ be a generic network causal effect. Consider the class of all linear estimators of $\theta$ under $P(\mathbf{Z})$. If $\theta$ is finite, then the Horvitz–Thompson estimator is inadmissible w.r.t. MSE in this class.*

The proof of Theorem 4.1 uses the next lemma: given any unbiased estimator whose minimal variance over possible $\theta$ is positive, a shrinkage estimator can be constructed with strictly smaller MSE uniformly.

**Lemma 4.1.** *Let $\mathbb{P} = P(\mathbf{Z})$ be any design and let $\hat{\theta}$ be any unbiased estimator of a finite causal effect $\theta$ under $\mathbb{P}$. If* $\min_\theta \mathrm{Var}_{\mathbb{P}}(\hat{\theta}) > 0$*, then there exists an estimator $\hat{\theta}_1$ with* $\mathrm{MSE}(\hat{\theta}_1) < \mathrm{MSE}(\hat{\theta})$ *for all $\theta$.*

Common designs such as Bernoulli and CRD are random exposure designs because they do not take into account the network structure, and hence have no direct control over the number of units exposed. These are the simplest designs that can be used in practice, as the interference structure is usually not known at the design stage. Hence the H–T estimator is inadmissible for average causal effects under these designs.

We now consider designs that fix the number of units at a given treatment–exposure pair. Such designs necessarily need to be aware of the interference graph, which is usually not known at the design stage, so they can be impractical to implement, or require the collection of additional information on the interference structure.

**Definition 3** (Fixed exposure designs)**.** *With notation as in Definition 2, a design $P(\mathbf{Z})$ is a* fixed exposure design *for $\theta$ if* $\mathrm{Var}(n_{d_1}) = 0$ *and* $\mathrm{Var}(n_{d_0}) = 0$.

**Theorem 4.2.** *For a fixed exposure design and any generic causal effect $\theta$, the H-T estimator is admissible (w.r.t. MSE) among all linear estimators.*

Fixed exposure designs can be implemented via restricted randomization schemes. One concrete example is the independent-set design of Karwa and Airoldi (2018) for estimating average direct effects under neighborhood interference. Construct an independent set of size $n_I$; randomly choose $n_1$ nodes from it and assign them to treatment while assigning their neighbors to control, with all remaining nodes assigned to control. This procedure fixes the counts for $(z_i = 1, e_i = 0)$ and $(z_i = 0, e_i = 0)$, yielding a fixed exposure design for the corresponding estimand.

**Remark 1.** *Theorem 4.1 and Theorem 4.2 have analogues in settings without network interference and in network sampling. For instance, in the classic no-interference setting, Bernoulli assignment with H–T estimator is an inadmissible estimation strategy (by a similar shrinkage argument), while CRD with the corresponding H–T estimator is part of an admissible strategy. In the network setting, however, the CRD paired with the H–T estimator is inadmissible because exposure counts are random, adding variance.*

## 5. Admissibility of the Hajek and Ratio estimators for random exposure designs

As shown in the previous section, the H–T estimator is inadmissible under random exposure designs but admissible under fixed exposure designs. Since fixed exposure designs are often impractical (they require knowing the interference graph), we ask whether better estimators exist under common random exposure designs. This section studies two such alternatives: the Hájek and the ratio estimators. We show that, under a mild condition on weights, both estimators are admissible (w.r.t. MSE) among all linear estimators for any design, including random exposure designs. In fact, we present a broad class of linear admissible estimators that contains both Hájek and ratio estimators.

**Lemma 5.1.** *Let $P(\mathbf{Z})$ be any design. Let $\hat{\theta} = \hat{\theta}(d_1) - \hat{\theta}(d_0)$ be a linear estimator of $\theta(d_1, d_0) = \theta(d_1) - \theta(d_0)$ with weights $w_i(\mathbf{Z}, d_j)$, where $\hat{\theta}(d_j) = (1/n)\sum_i Y_i(d_j)w_i(\mathbf{Z}, d_j)$. If,*

1. *For each $j \in \{0, 1\}$, there exists a set of strictly positive potential outcomes $\{Y_i^*(d_j)\}$, such that $MSE(\hat{\theta}(d_j), \theta(d_j)) = 0$ at $\{Y_i^*(d_j)\}$,*
2. *For each $i$, $j \in \{0, 1\}$, $\min_{\mathbf{z} \in \Omega_i(d_j)} w_i(\mathbf{z}, d_j) > 1$ where $\Omega_i(d_j) = \{\mathbf{z} : I_i(d_j) = 1\}$*
3. *For each $i$, $j \in \{0, 1\}$, $\pi_i(d_j) = P(\mathbf{z} \in \Omega_i(d_j)) > 0$,*

*then, $\hat{\theta}$ is admissible with respect to mean-square error in the class of all linear weighted estimators, where the weights are allowed to depend on the entire vector* $\mathbf{z}$.

Lemma 5.1 gives a wide class of admissible linear estimators: if an estimator attains zero MSE at some positive potential-outcome vector and its weights on the set of treatment allocations that reveal the potential outcome for unit $i$ are greater than 1 (and propensities are positive), then it is admissible. This is the contrapositive of Lemma 4.1, which constructs a shrinkage estimator whenever an estimator has strictly positive variance.

**Remark 2.** *Lemma 5.1 requires* $\min_{\mathbf{z}\in\Omega_i(d_j)} w_i(\mathbf{z}, d_j) > 1$ *for each* $i, j$. *A practical way to ensure this is via restricted randomization tailored to the estimator; alternatively, one can view admissibility conditionally on the event where the weights exceed one.*

### 5.1. The Hájek estimator

Define the H–T estimates of counts and sums:

$$\hat{n}(d_1) = \sum_{i=1}^{n} \frac{I(D_i = d_1)}{\pi_i(d_1)}, \qquad \hat{n}(d_0) = \sum_{i=1}^{n} \frac{I(D_i = d_0)}{\pi_i(d_0)},$$

$$\hat{S}(d_1) = \sum_{i=1}^{n} \frac{I(D_i = d_1)}{\pi_i(d_1)} Y_i(D_i), \qquad \hat{S}(d_0) = \sum_{i=1}^{n} \frac{I(D_i = d_0)}{\pi_i(d_0)} Y_i(D_i).$$

The Hájek estimator (Hájek, 1971) is

$$\hat{\theta}_J = \frac{\hat{S}(d_1)}{\hat{n}(d_1)} - \frac{\hat{S}(d_0)}{\hat{n}(d_0)}.$$

It is defined when $\hat{n}(d_j) > 0$. The Hájek estimator is biased but asymptotically unbiased; see Aronow and Samii (2017). Under the sufficient condition $\hat{n}(d_j) < n$ for $j = 0, 1$ (for all relevant allocations), the Hájek estimator satisfies the weight condition in Lemma 5.1 and is therefore admissible:

**Theorem 5.1.** *Consider a design such that (a)* $\pi_i(d_j) > 0$ *for all* $i, j \in \{0, 1\}$, *and (b) for all* $\mathbf{z} \in \Omega$, $\hat{n}(d_1) < n$ *and* $\hat{n}(d_0) < n$. *Then* $\hat{\theta}_J$ *is admissible (w.r.t. MSE) among all linear estimators for* $\theta$.

The condition $\hat{n}(d_j) < n$ is sufficient but not necessary; admissibility may hold in broader cases. Interpreting Theorem 5.1 conditionally (on the event $E = \{\mathbf{Z} : \hat{n}(d_j) < n\}$) or using restricted randomization (e.g., re-randomization, see Morgan and Rubin (2012)) can make this condition practical. Theorem 5.2 below gives an analytic approximation showing $P(\hat{n}(d_j) < n)$ can be appreciable when propensities $\pi_i$ are less than $1/2$. Even though $E[\hat{n}(d_j)] = n$, the distribution of $\hat{n}(d_j)$ can be skewed and the skewness depends on the propensity scores $\pi_i(d_j)$ - for smaller propensity scores the distribution is right skewed making the $P(\hat{n}(d_j) < n) > 1/2$.

**Theorem 5.2.** *Consider a Bernoulli design. For fixed* $d_j$, *let* $I_i = \mathbf{1}\{D_i = d_j\}$ *and* $\pi_i = P(I_i = 1)$. *Let* $G$ *be the interference graph and denote* $d_i = |N_i|$. *Assume (a)*

*$d_i$ are uniformly bounded, (b) there exist constants $c$ and $C$ with $c < \pi_i < C$, and (c)* $\mathrm{Var}(\hat{n}) \geq cn$ *for some $c > 0$. Then, as $n \to \infty$,*

$$P\big(\widehat{n}(d_j) < n\big) \;\approx\; \frac{1}{2} \;+\; \frac{\gamma^*}{6\sqrt{2\pi}},$$

*where $\gamma^* = (\sigma^*)^{-3/2} \sum_i \frac{(1-\pi_i)(1-2\pi_i)}{\pi_i^2}$ is the standardized skewness of $\sum_i I_i^*/\pi_i$ where $I_i^* \sim I_i$ are independent. In particular, if $\pi_i < 1/2$, then $P(\hat{n}(d_j) < n) \gtrapprox 1/2$.*

The proof approximates the dependent weighted sum $\hat{n} = \sum_i I_i/\pi_i$ by a independent weighted sum $\sum_i I_i^*/\pi_i$ via Chen and Shao (2004) and then applies a one-term Edgeworth expansion to the independent sum using a weak Cramér condition of Angst and Poly (2017).

**Remark 3.** *Theorem 5.2 uses a normal/Edgeworth approximation but may not cover all practical cases, for e.g., highly heterogeneous or very small propensities; nevertheless $P(\hat{n}(d) < n)$ can still be large in such settings.*

### *5.2. The ratio estimator*

Let $n_{d_1} = \sum_i I(D_i = d_1)$ and $n_{d_0} = \sum_i I(D_i = d_0)$ be observed counts. The ratio estimator is

$$\hat{\theta}_r = \frac{1}{n}\left(\frac{\sum_i \pi_i(d_1)}{n_{d_1}}\hat{S}(d_1) - \frac{\sum_i \pi_i(d_0)}{n_{d_0}}\hat{S}(d_0)\right).$$

Note $n_{d_j}$ is different from $\hat{n}(d_j)$. In fixed exposure designs $n_{d_j} = E[n_{d_j}] = \sum_i \pi_i(d_j)$, so the ratio estimator equals H–T estimator. Thus the ratio estimator generalizes the H–T estimator to random exposure designs and is self-normalizing like the Hájek estimator.

**Theorem 5.3.** *If a design yields $n_{d_1} < E[n_{d_1}]$ and $n_{d_0} < E[n_{d_0}]$, then $\hat{\theta}_r$ is admissible (w.r.t. MSE) among all linear estimators for $\theta$.*

As with Hájek, Theorem 5.3 is naturally interpreted conditionally and analytic approximations for $P(n_{d_j} < E[n_{d_j}])$ can be derived similar to Theorem 5.2.

We illustrate these results numerically using the canonical $2 \times 2$ exposure model with a Bernoulli design. Each unit is exposed iff at least one neighbor is treated, so $d \in \{(1,1), (1,0), (0,1), (0,0)\}$. Potential outcomes follow

$$Y_i(z_i, e_i) = \alpha_i + \beta_i z_i + \gamma_i e_i + \delta_i z_i e_i,$$

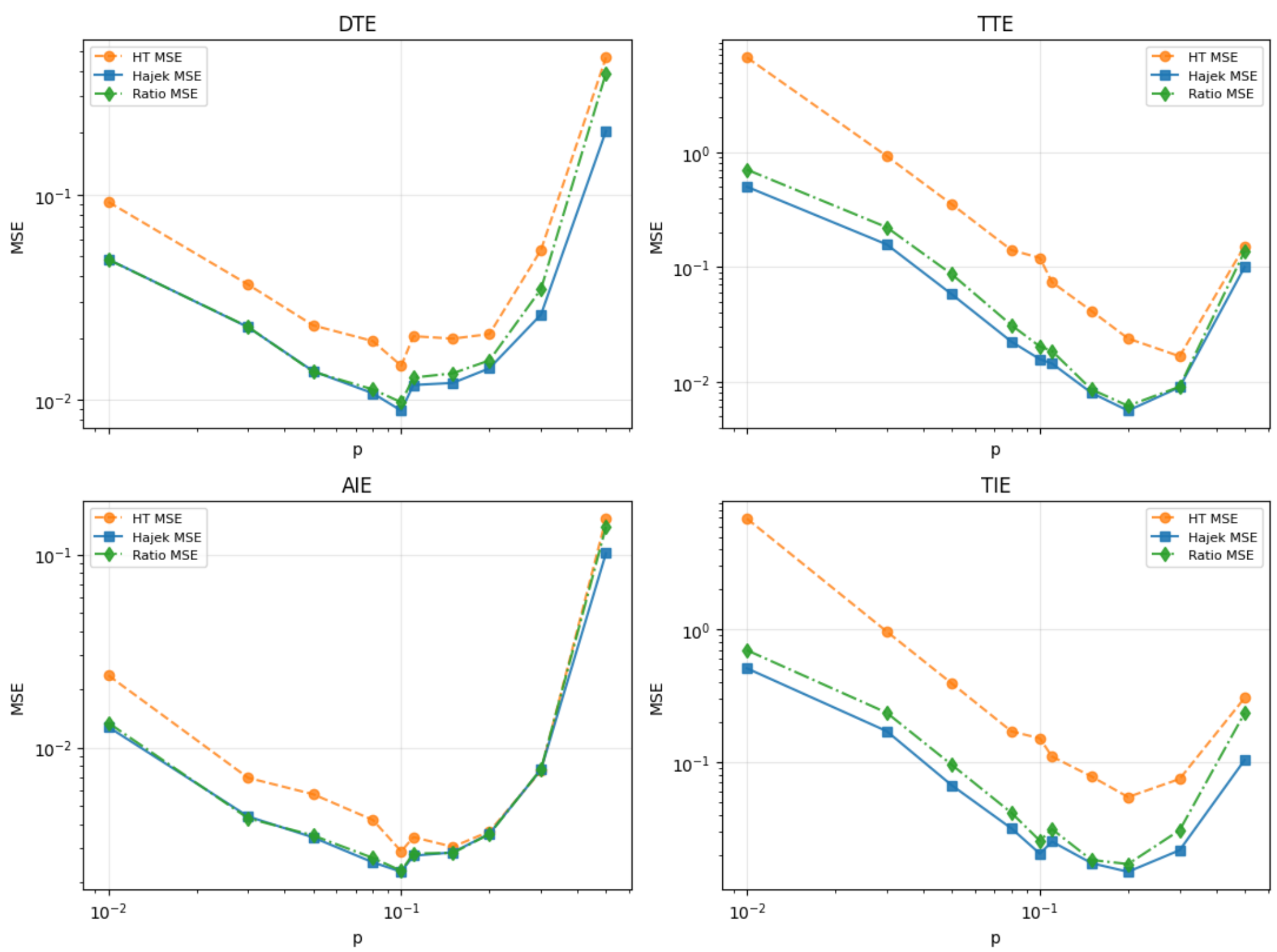

*Fig 1. MSE of HT, Hajek and Ratio estimators under the* 2 *by* 2 *exposure model for 4 different estimands: DTE, TTE, AIE and TIE. The design is a Bernoulli design with different values of* $p = \{0.01, 0.03, 0.05, 0.08, 0.1, 0.15, 0.20, 0.30, 0.50\}$. *The number of units is* $n = 2000$.

with independent coefficients $\alpha_i \sim \mathcal{N}(0, 1)$, $\beta_i \sim \mathcal{N}(1, 0.2^2)$, $\gamma_i \sim \mathcal{N}(0.5, 0.2^2)$, $\delta_i \sim \mathcal{N}(0.2, 0.1^2)$. We use $n = 2000$, Bernoulli $p \in \{0.01, 0.03, 0.05, 0.08, 0.1, 0.15, 0.2, 0.3, 0.5\}$, and a bounded-degree interference graph (degrees uniform in $\{2, \ldots, 9\}$) generated via the configuration model. For each $p$ we perform $R = 200$ randomizations and compute HT, Hájek, and ratio estimators for four estimands (1) Direct Treatment Effect (DTE) = $\tau(1, 0) - \tau(0, 0)$, (2) Total treatment effect (TTE) = $\tau(1, 1) - \tau(0, 0)$, (3) Average interference effect (AIE) = $\tau(0, 1) - \tau(0, 0)$, (4) Total interference effect (TIE) = $\tau(1, 1) - \tau(1, 0)$, recording the MSEs. Across simulations the Hájek estimator attains the lowest MSE, the ratio estimator is next, and $P(\hat{n}(d_j) < n)$ and $P(n_{d_j} < E[n_{d_j}])$ ranged roughly 0.45–0.6. Figure 1 shows the results.

## 6. A new admissible estimator for random exposure designs

While the Hájek and ratio estimators are admissible under random exposure designs, they can be biased. We now introduce a new estimator that is both unbiased and admissible under *any* design, including random exposure designs. The key idea is to analyze the design conditional on an ancillary statistic—the random number of units assigned to each treatment–exposure combination.

To illustrate, consider the no-interference case with a Bernoulli design for estimating the average treatment effect. Here the numbers of treated and control units are random, so by the results of Section 4, the H–T estimator under the Bernoulli design is inadmissible. However, conditioning the propensity scores on the realized counts of treated and control units effectively converts the design to one with fixed sample sizes. The resulting conditional estimator—equivalent to the difference-in-means estimator—is admissible, whereas the unconditional H–T estimator is not (see Remark 4). This conditionalization of propensity scores extends naturally to settings with interference.

*The conditional Horvitz–Thompson estimator:* Let $n_{d_1} = \sum_i I(D_i = d_1)$ and $n_{d_0} = \sum_i I(D_i = d_0)$ denote the (possibly random) numbers of units observed under $d_1$ and $d_0$. Define the *conditional propensity scores*

$$\pi_i(D_i = d_1 \mid n_{d_1}) = \sum_{\mathbf{z} \in \Omega_i(d_1 | n_{d_1})} P(\mathbf{Z}), \qquad \pi_i(D_i = d_0 \mid n_{d_0}) = \sum_{\mathbf{z} \in \Omega_i(d_0 | n_{d_0})} P(\mathbf{Z}),$$

where $\Omega_i(d_1 \mid n_{d_1} = t) = \{\mathbf{z} : \sum_j I(z_j = z, e_j = e) = t,\ z_i = z, e_i = e\}$ and $\Omega_i(d_0 \mid n_{d_0} = t) = \{\mathbf{z} : \sum_j I(z_j = z', e_j = e') = t,\ z_i = z', e_i = e'\}$. Each $\pi_i(D_i = d_j \mid n_{d_j})$ is random and it depends on the full assignment vector $\mathbf{Z}$. The *conditional H–T estimator* uses these conditional propensities and is defined as:

$$\hat{\theta}_{CHT} = \frac{1}{n}\left(\sum_{i=1}^{n} \frac{I(D_i = d_1)}{\pi_i(D_i = d_1 \mid n_{d_1})} Y_i(D_i) - \sum_{i=1}^{n} \frac{I(D_i = d_0)}{\pi_i(D_i = d_0 \mid n_{d_0})} Y_i(D_i)\right).$$

Note that for fixed exposure designs, $\hat{\theta}_{CHT}$ coincides with the usual H–T estimator.

**Proposition 6.1.** *If $\pi_i(d_j \mid n_{d_j}) > 0$ for all $i, j$, then $\mathbb{E}[\hat{\theta}_{CHT}] = \theta(d_1, d_0)$, where the expectation is under $P(\mathbf{Z})$.*

**Theorem 6.1.** *If $0 < \pi_i(d_j \mid n_{d_j}) < 1$ for all $i$ and $j \in \{0, 1\}$, then the conditional H–T estimator is admissible (w.r.t. MSE) in the class of all linear estimators under any design.*

Together, Proposition 6.1 and Theorem 6.1 imply that whenever the conditional propensity scores are bounded away from 0 and 1, the conditional H–T estimator is both unbiased and admissible for any design.

**Remark 4.** *Under SUTVA, potential outcome of unit $i$ depends only on $z_i$. In a Bernoulli design with treatment probability $p$, the numbers of treated and control units are random, making it a random sample size design; a completely randomized design (CRD) with fixed treated count $n_t$ is a fixed sample size design. For the Bernoulli design, the conditional propensities are $(\sum_i Z_i)/n$ and $(\sum_i(1 - Z_i))/n$ for treated and control units, versus the usual $p$ and $1 - p$. The conditional H–T estimator becomes*

$$\hat{\theta}_{CHT}(B) = \frac{\sum_i y_i(Z_i)Z_i}{\sum_i Z_i} - \frac{\sum_i y_i(Z_i)(1 - Z_i)}{\sum_i(1 - Z_i)},$$

*while the ordinary H–T estimator is*

$$\hat{\theta}_{HT}(B) = \frac{1}{n}\Big(\sum_i \frac{y_i(Z_i)Z_i}{p} - \sum_i \frac{y_i(Z_i)(1 - Z_i)}{1 - p}\Big).$$

*Thus, under Bernoulli designs (which is a random sample size design), $\hat{\theta}_{CHT}(B)$ is admissible (and reduces to the difference-in-means estimator), whereas $\hat{\theta}_{HT}(B)$ is not admissible. On the other hand, under a CRD (which is a fixed sample size design), both estimators coincide and are admissible:*

$$\hat{\theta}_{HT}(CRD) = \hat{\theta}_{CHT}(CRD) = \frac{\sum_i y_i(Z_i)Z_i}{n_t} - \frac{\sum_i y_i(Z_i)(1 - Z_i)}{n_c}.$$

### *6.1. Approximating the conditional H–T estimator*

Under network interference the conditional H–T estimator is harder to compute than in SUTVA because the conditional propensity scores $\pi_i(d_j \mid n_{d_j})$ involve dependent indicators. We present three approximations: (i) a Poisson (or Binomial) analytic approximation under a Bernoulli design, (ii) a simple analytic approximation when the conditional propensities are nearly equal or there is exchangeability, and (iii) a numerical MCMC approximation.

Fix $j$ and write $I_i = \mathbf{1}\{D_i = d_j\}$ and $\underline{I} = (I_1, \ldots, I_n)$ with marginals $\pi_i = P(I_i = 1)$. Under a Bernoulli design $I_i$ and $I_k$ are dependent only when nodes $i$ and $k$ lie within graph distance $\leq 2$. Let $Z_{N_i}(2) = \{k : d(i, k) \leq 2\}$ denote the 2-hop neighborhood of $i$. Then $\sum_i I_i$ is a sum of dependent Bernoulli variables (a dependent Poisson–Binomial). When the 2-hop neighborhoods are small (grow as $o(n)$) we can approximate $\sum_i I_i$ by a Poisson or Binomial law and compute conditional probabilities $P(I_i = 1 \mid \sum_i I_i = t)$ analytically. If the conditional propensities themselves are nearly equal, or the interference graph and the exposure functions are exchangeable, a second analytic approximation gives a simple reduction to difference-in-means. Finally, one can approximate the dependent vector by independent Bernoullis and use MCMC to estimate conditional probabilities. The next lemma formalizes the first approximation.

**Lemma 6.1.** *Consider a Bernoulli design. Let* $Z_{N_i}(2) = \{j \in [n] : d(j,i) \leq 2\}$. *If (a)* $\max_i |Z_{N_i}(2)|/n \to 0$ *and (b)* $\max_i n\pi_i \to c \in (0,\infty)$, *then there exist independent Poisson random variables* $W_1, \ldots, W_n$ *with* $E[W_i] = P(I_i = 1)$ *such that*

$$\left|P(I_i = 1 \mid \sum_i I_i) - P(W_i = 1 \mid \sum_i W_i)\right| \to 0.$$

Lemma 6.1 yields a simple analytic approximation to the conditional H–T estimator: when the Poisson approximation applies (and $\pi_i/\sum_j \pi_j \to 0$), the conditional estimator is approximated by the familiar ratio estimator, as shown by the corollary below.

**Corollary 1.** *Under Lemma 6.1's conditions and* $\pi_i/\sum_j \pi_j \to 0$,

$$\hat{\theta}_{CHT} \approx \frac{E[n_{d_1}]}{n_{d_1}} \hat{\theta}_{HT}(d_1) - \frac{E[n_{d_0}]}{n_{d_0}} \hat{\theta}_{HT}(d_0),$$

*where* $\hat{\theta}_{HT}(d_j) = \frac{1}{n} \sum_{i=1}^{n} \frac{Y_i(D_i)I(D_i=d_j)}{\pi(D_i=d_j)}$.

The second method approximates the CHT by a difference-in-means estimator, when the conditional propensity scores are approximately equal. This can happen, for example, under exchangeable designs (Bernoulli or CRD), when the exposure function depend only on the number of treated neighbors, and interference neighborhoods are similar in size and structure.

**Proposition 6.2.** *Let* $p_i(j) = \pi_i(d_j \mid n_{d_j})$ *and* $\bar{p}(j) = \frac{1}{n}\sum_i p_i(j)$. *Define* $\beta_n(j) = \frac{1}{n}\sum_{i=1}^{n} |p_i(j) - \bar{p}(j)|$, $\quad \beta_n = \max_j \beta_n(j)$. *If potential outcomes are uniformly bounded,* $\min_i p_i(j) > 0$, *and* $\bar{p}(j) > c > 0$, *then* $|\hat{\theta}_{CHT} - \hat{\theta}_{df}| = O(\beta_n)$, *where*

$$\hat{\theta}_{df} = \frac{\sum_i Y_i(D_i)I(D_i = d_1)}{\sum_i I(D_i = d_1)} - \frac{\sum_i Y_i(D_i)I(D_i = d_0)}{\sum_i I(D_i = d_0)}$$

*is the difference-in-means estimator.*

Finally, we describe an MCMC approximation that works well in practice. Lemma 6.1 suggests approximating $(I_1, \ldots, I_n)$ by independent Poissons, but it turns out that the independent Bernoulli vector $(V_1, \ldots, V_n)$ with $P(V_i = 1) = \pi_i$ often performs better. To estimate $P(V_i = 1 \mid \sum_i V_i = t)$ we sample from the conditional distribution of $(V_1, \ldots, V_n) \mid \sum_i V_i$ using a simple switch Markov chain of Heng et al. (2020), which mixes in $O(n \log n)$ time.

## 7. Numerical Results

We illustrate our theoretical findings on a real-world dataset. We compare the performance of the following six estimators: 1) H-T estimator, 2) Hajek estimator 3) Ratio estimator 4) Difference-in-means estimator 5) Approximate Conditional H-T estimator that uses the MCMC algorithm to estimate the conditional propensity scores and (6) Approximate Conditional Hajek estimator. The approximate conditional Hajek estimator is defined as the Hajek estimator that uses the approximate conditional propensity scores. We use friendship networks from the National Longitudinal Study of Adolescent Health ("Add Health") studied in Aronow and Samii (2017). The dataset contains 144 networks (average size ~600). We symmetrize edges (an undirected edge exists if either student nominated the other) and drop isolated nodes; see Aronow and Samii (2017) for details.

We adopt a binary $2 \times 2$ treatment–exposure model where a unit is exposed if at least one neighbor is treated ($E_i = \mathbf{1}\{\sum_{j \in N_i} z_j > 1\}$). Potential outcomes follow the nonparametric form $Y_i(z_i, e_i) = \alpha_i + \beta_i z_i + \gamma_i e_i + \delta_i z_i e_i$, so each unit has four potential outcomes $y_i(1,1), y_i(1,0), y_i(0,1), y_i(0,0)$. We treat the observed number of after-school activities as the baseline $y_i(0,0)$ and use the "dilated effects" specification of Rosenbaum (1999); Aronow and Samii (2017) to set $y_i(1,1) = 2y_i(0,0), \quad y_i(1,0) = 1.5y_i(0,0), \quad y_i(0,1) = y_i(0,0)$. We estimate four causal contrasts: Total treatment effect = $\bar{y}(1,1) - \bar{y}(0,0) = \bar{\beta} + \bar{\gamma} + \bar{\delta}$, Average direct effect $\bar{y}(1,0) - \bar{y}(0,0) = \bar{\beta}$, Additive interference $\bar{y}(0,1) - \bar{y}(0,0) = \bar{\gamma}$, and Total interference $\bar{y}(1,1) - \bar{y}(1,0) = \bar{\gamma} + \bar{\delta}$.

For each network we draw $R = 500$ Bernoulli randomizations with treatment probability $p = 0.01$. For each randomization we compute all six estimators for the four estimands and then aggregate results across networks. Table 1 reports bias, standard deviation (sd), root mean squared error (RMSE), and the Monte Carlo standard error of the bias (bstderr). The results are shown in Table 1. The results suggest the following: (1) For all estimands (except $\bar{\gamma}$), the approximate-conditional H-T estimator always has considerably lower RMSE than the H-T estimator. For $\bar{\gamma}$, the two estimators are very close. (2) For all the estimands, the Hajek estimator and the approximate-conditional Hajek estimator have the least RMSE, and are comparable to each other. (3) For estimating the total treatment effect and the average direct effect, the ratio estimator, Hajek estimator and the approximate-conditional Hajek estimator have similar behavior and have a smaller RMSE. (4) For estimating the interference effects, the conditional Hajek and the Hajek estimators have similar performance, both giving the least possible RMSE.

*Table 1*
*Results of the simulation study on Add health networks*

| Estimand | Estimator | Bias | sd | RMSE | bstderr |
|---|---|---|---|---|---|
| Total treatment effect ($\bar{\beta}+\bar{\gamma}+\bar{\delta}$) | a-conditional HT | -0.38 | 0.74 | 0.53 | 0.03 |
| | a-conditional Hajek | -0.04 | 0.71 | 0.16 | 0.01 |
| | Horvitz Thompson | -0.45 | 1.65 | 1.47 | 0.12 |
| | Difference-in-means | -0.49 | 0.77 | 0.55 | 0.02 |
| | Hajek | -0.07 | 0.7 | 0.16 | 0.01 |
| | Ratio | -0.04 | 0.73 | 0.16 | 0.01 |
| Average direct effect ($\bar{\beta}$) | a-conditional HT | -0.27 | 0.4 | 0.35 | 0.02 |
| | a-conditional Hajek | 0.02 | 0.33 | 0.06 | 0.01 |
| | Horvitz Thompson | -0.06 | 0.52 | 0.36 | 0.03 |
| | Difference-in-means | 0.07 | 0.33 | 0.1 | 0.01 |
| | Hajek | 0.01 | 0.33 | 0.06 | 0.01 |
| | Ratio | -0.03 | 0.37 | 0.09 | 0.01 |
| Additive interference effect ($\bar{\gamma}$) | a-conditional HT | -0.09 | 0.24 | 0.18 | 0.01 |
| | a-conditional Hajek | 0.02 | 0.2 | 0.07 | 0.01 |
| | Horvitz Thompson | -0.01 | 0.21 | 0.12 | 0.01 |
| | Difference-in-means | -0.37 | 0.29 | 0.41 | 0.02 |
| | Hajek | -0.02 | 0.18 | 0.05 | 0 |
| | Ratio | -23.92 | 6.96 | 24.87 | 0.58 |
| Total Interference effect ($\bar{\gamma}+\bar{\delta}$) | a-conditional HT | -0.12 | 0.42 | 0.3 | 0.02 |
| | a-conditional Hajek | -0.07 | 0.42 | 0.21 | 0.02 |
| | Horvitz Thompson | -0.4 | 1.34 | 1.31 | 0.11 |
| | Difference-in-means | -0.57 | 0.5 | 0.64 | 0.03 |
| | Hajek | -0.09 | 0.41 | 0.21 | 0.02 |
| | Ratio | -5.99 | 3.79 | 6.95 | 0.3 |

## 8. Implications for Practitioners

*Estimating network causal effects:* The H–T estimator is the standard tool for estimating causal effects under network interference but, as shown in this paper, it is inadmissible with respect to MSE under random exposure designs. In practice, such designs are the most common because the number of units assigned to each treatment–exposure combination is typically not directly controlled. Practitioners therefore face two options: (i) implement a fixed exposure design and use the H–T estimator, achieving unbiasedness and admissibility, or (ii) use a random exposure design but replace the H–T estimator with an admissible alternative such as the Hájek, ratio, or conditional H–T estimator (or its approximations). These estimators are especially suitable for estimating network causal effects like the total treatment, total interference, and average direct effects.

*Implications for causal inference without interference:* The results also apply when interference is absent. Theorem 4.1 implies that the H–T estimator under a Bernoulli design is inadmissible, whereas it becomes admissible under a completely randomized design (CRD). Under Bernoulli designs, the conditional H–T estimator provides an admissible alternative. Because many observational studies are analyzed as Bernoulli designs—using estimated propensity scores as inverse weights—our results suggest that conditional H–T, ratio, or Hájek estimators may be preferable to the classical H–T estimator in these cases. The inefficiency of the H–T arises from extra variance due to the random numbers of treated and control units, a problem mitigated by conditional or ratio estimators. This phenomenon parallels the finding of Robins et al. (1994), where using estimated rather than true propensity scores improved efficiency. As noted in Remark 4, conditional propensity scores in the Bernoulli design play a similar role as estimated propensity scores. The same reasoning extends to network sampling methods such as respondent-driven or snowball sampling, where the number of sampled units is random. In these cases, the H–T estimator is likewise inadmissible with respect to MSE.

## 9. Discussion

We study the Horvitz–Thompson (H–T) estimator for estimating causal effects under network interference by embedding it in a broad class of linear weighted estimators and analyzing optimality properties. We characterize the set of all possible linear unbiased estimators and show that there can be infinitely many such estimators, depending on the exposure model. when weights satisfy a natural sufficiency constraint, the H–T estimator is the unique unbiased estimator. However, unbiasedness relies on correct specification of the exposure model and knowledge of the interference graph—both often unavailable at design time. Within the larger class of linear estimators, H–T is inadmissible when the design induces ran-dom counts for the treatment–exposure combination of interest (so-called *random exposure* designs) and admissible when those counts are fixed ("*fixed exposure* designs"). Fixed exposure designs remove the extra randomness by construction but typically require knowledge of the network; we gave examples based on restricted randomization.

We also show that under a mild condition, the Hajek estimator, and the ratio estimator are admissible under any design. We do this by characterizing a class of admissible estimators. This class relies on the key fact that if the minimum variance of the estimator is strictly positive, then it is inadmissible. On the other hand, if there exists a set of potential outcomes where the variance of the estimator is 0, then it is admissible. Using this insight, we propose a new estimator called the conditional Horvitz Thompson estimator that uses propensity scores conditioned on the number of units under a fixed treatment-exposure combination. This conditioning allows the minimum variance of the H-T estimator to be zero, even for a random exposure design, thus making it admissible.

There are many key open questions that remain. A central open question is to study admissibility with respect to the class of all estimators, and to find admissible estimators under random exposure designs without imposing the constraint of weights being greater than 1. A related question is to construct fixed exposure designs, such designs necessarily depend on the network structure and the estimand. Another direction of problems is to find conditions under which average treatment effects can be estimated when the interference graph and the exposure model are not known. Some progress has been made in this direction, see Yu et al. (2022); Sävje (2023). One possibility is to consider estimators and designs that are robust to the interference graph and the exposure model, another would be to learn the interference graph and the exposure model from the data. While there is some work on testing for the existence of interference (Aronow, 2012; Pouget-Abadie et al., 2019), an important related question that deserves further investigation is testing the assumed form of interference. Finally, in this paper, we considered the problem of admissibility for a fixed design. An important direction for future work is to consider estimator and design pairs that are *uniformly* admissible.

## Appendix A: Weights of the Horvitz-Thompson Estimator for common exposure models

The weight of an H-T estimator for unit $i$ is inversely proportional to the propensity score $\pi_i(z, e)$, the probability of observing that potential outcome under the design $P(\mathbf{Z})$. These probabilities depend on the design and the exposure model. In this section, we compute analytical formulae of these probabilities for the CRD and the Bernoulli designs for different exposure models.

**Theorem A.1** (Propensity Scores for Symmetric Exposure)**.** *Consider the symmetric exposure function,* $e_i = f(\mathbf{Z}_{N_i}) = |\mathbf{Z}_{N_i}|$*,* $e_i \in \{0, 1, \ldots, d_i\}$*. For a CRD Design,*

$$\mathbb{P}\left(Z_i = 1, E_i = e_i\right) = \frac{n_t}{n}\frac{\binom{n_t-1}{e_i}\binom{n_c}{d_i-e_i}}{\binom{n-1}{d_i}} \text{ if } n_t \geq e_i + 1 \text{ and } n_c \geq d_i - e_i, 0 \text{ otherwise}$$

$$\mathbb{P}\left(Z_i = 0, E_i = e_i\right) = \frac{n_c}{n}\frac{\binom{n_t}{e_i}\binom{n_c-1}{d_i-e_i}}{\binom{n-1}{d_i}} \text{ if } n_t \geq e - i \text{ and } n_c \geq d_i - e_i + 1, 0 \text{ otherwise}$$

*For a Bernoulli Design,*

$$\mathbb{P}\left(Z_i = 1, E_i = e_i\right) = \binom{d_i}{e_i} p^{e_i+1}(1-p)^{d_i-e_i}$$

$$\mathbb{P}\left(Z_i = 0, E_i = e_i\right) = \binom{d_i}{e_i} p^{e_i}(1-p)^{d_i-e_i+1}$$

**Theorem A.2** (Propensity Scores for Binary Exposure)**.** *Consider the symmetric exposure function,* $e_i = f(\mathbf{Z}_{N_i}) = I(|\mathbf{Z}_{N_i}| \geq 1)$*,* $e_i \in \{0, 1\}$*. i.e a unit is exposed if* *at least 1 of its neighbor is treated. For a CRD,*

$$\mathbb{P}\left(Z_i = 1, E_i = 1\right) = \begin{cases} 0 & \text{if } d_i = 0 \\ \frac{n_t}{n}\left[1 - \frac{\binom{n_c}{d_i}}{\binom{n-1}{d_i}}\right] & \text{if } 0 < d_i \leq n_c \\ \frac{n_t}{n}, & \text{if } d_i > n_c \end{cases}$$

$$\mathbb{P}\left(Z_i = 1, E_i = 0\right) = \begin{cases} \frac{n_t}{n} & \text{if } d_i = 0 \\ \frac{n_t}{n}\frac{\binom{n_c}{d_i}}{\binom{n-1}{d_i}} & \text{if } 0 < d_i \leq n_c \\ 0, & \text{if } d_i > n_c \end{cases}$$

$$\mathbb{P}\left(Z_i = 0, E_i = 1\right) = \begin{cases} 0 & \text{if } d_i = 0 \\ \frac{n_c}{n}\left[1 - \frac{\binom{n_c-1}{d_i}}{\binom{n-1}{d_i}}\right] & \text{if } 0 < d_i \leq n_c - 1 \\ \frac{n_c}{n}, & \text{if } d_i > n_c - 1 \end{cases}$$

$$\mathbb{P}\left(Z_i = 0, E_i = 0\right) = \begin{cases} \frac{n_c}{n} & \text{if } d_i = 0 \\ \frac{n_c}{n}\frac{\binom{n_c-1}{d_i}}{\binom{n-1}{d_i}} & \text{if } 0 < d_i \leq n_c - 1 \\ 0 & \text{if } d_i > n_c - 1 \end{cases}$$

*Similarly, for a Bernoulli trial with probability of success p, we have*

$$\mathbb{P}\left(Z_i = 1, E_i = 1\right) = p(1 - (1-p)^{d_i})$$
$$\mathbb{P}\left(Z_i = 1, E_i = 0\right) = p(1-p)^{d_i}$$
$$\mathbb{P}\left(Z_i = 0, E_i = 1\right) = (1-p)(1 - (1-p)^{d_i})$$
$$\mathbb{P}\left(Z_i = 0, E_i = 0\right) = (1-p)^{d_i+1}$$

*Under a cluster randomized design, let $u_i$ be the number of clusters of unit i and*

*it's neighbors. Assume* $n_k > 0, \forall k = 1, \ldots, K$.

$$
\begin{aligned}
&\mathbb{P}\left(z_i = 1, e_i = 1\right) = \frac{K_t}{K} \\
&\mathbb{P}\left(z_i = 1, e_i = 0\right) = 0 \\
&\mathbb{P}\left(z_i = 0, e_i = 1\right) = 0 \textit{ if } u_i = 1, \frac{K_c}{K}\left[1 - \prod_{i=1}^{u_i-1} \frac{K_c - u_i}{K - u_i}\right] \textit{ if } u_i > 1 \\
&\mathbb{P}\left(z_i = 0, e_i = 0\right) = \frac{K_c}{K} \textit{ if } u_i = 1, \frac{K_c}{K}\left[\prod_{i=1}^{u_i-1} \frac{K_c - i}{K - i}\right] \textit{ if } u_i > 1 \\
&\qquad\qquad\qquad\quad = \prod_{i=1}^{u_i} \frac{K_c - i + 1}{K - i + 1}
\end{aligned}
$$

The results of this section demonstrate that, apart from the design, the weights of the H-T estimator also depend on the exposure model and the interference graph $G$. In particular, the unbiasedness property of the H-T estimator is sensitive to the mis-specification of interference graph $G$ and the exposure model $f(Z_{N_i})$. In many applications, $G$ and $f(Z_{N_i})$ are not known at the design stage.

## Appendix B: Missing proofs from Section 3

### B.1. Proof of Theorem 3.1

*Proof.* Note that by the consistency assumption, we have,

$$
\hat{\theta}_1 = \sum_{i=1}^{n} Y_i^{obs} w_i(\mathbf{z}) = \sum_{i=1}^{n} \sum_{z^*, e^*} Y_i(z^*, e^*) I(z_i = z^*, e_i = e^*) w_i(\mathbf{z})
$$

$$\begin{aligned}
\mathbb{E}[\hat{\theta}_1] &= \sum_{i=1}^{n} \sum_{z^*,e^*} Y_i(z^*,e^*) \left( \sum_{\mathbf{z}\in\Omega} I(z_i = z^*, e_i = e^*) w_i(\mathbf{z}) P(\mathbf{Z}) \right) \\
&= \sum_{i=1}^{n} \sum_{\mathbf{z}\in\Omega} w_i(\mathbf{z}) Y_i(z,e) I(z_i = z, e_i = e) P(\mathbf{Z}) \\
&\quad + \sum_{i=1}^{n} \sum_{\mathbf{z}\in\Omega} w_i(\mathbf{z}) Y_i(z',e') I(z_i = z', e_i = e') P(\mathbf{Z}) \\
&\quad + \sum_{i=1}^{n} \sum_{(z^*,e^*)\neq(z,e),(z',e')} \sum_{\mathbf{z}\in\Omega} w_i(\mathbf{z}) Y_i(z^*,e^*) I(z_i = z^*, e_i = e^*) P(\mathbf{Z}) \\
&= \sum_{i=1}^{n} \sum_{\mathbf{z}\in\Omega_i(z,e)} w_i(\mathbf{z}) Y_i(z,e) P(\mathbf{Z}) + \sum_{i=1}^{n} \sum_{\mathbf{z}\in\Omega_i(z',e')} w_i(\mathbf{z}) Y_i(z',e') P(\mathbf{Z}) \\
&\quad + \sum_{i=1}^{n} \sum_{(z^*,e^*)\neq(z,e),(z',e')} \sum_{\mathbf{z}\in\Omega_i(z^*,e^*)} w_i(\mathbf{z}) Y_i(z,e) P(\mathbf{Z}) \\
&= \sum_{i=1}^{n} Y_i(z,e) \left( \sum_{\mathbf{z}\in\Omega_i(z,e)} w_i(\mathbf{z}) P(\mathbf{Z}) \right) + \sum_{i=1}^{n} Y_i(z',e') \left( \sum_{\mathbf{z}\in\Omega_i(z',e')} w_i(\mathbf{z}) P(\mathbf{Z}) \right) \\
&\quad + \sum_{i=1}^{n} \sum_{(z^*,e^*)\neq(z,e),(z',e')} Y_i(z,e) \left( \sum_{\mathbf{z}\in\Omega_i(z^*,e^*)} w_i(\mathbf{z}) P(\mathbf{Z}) \right) \\
&= \sum_{i=1}^{n} Y_i(z,e) \left( \frac{1}{n} \right) - \sum_{i=1}^{n} Y_i(z',e') \left( \frac{1}{n} \right)
\end{aligned}$$

where the last line is required for unbiasedness. Since this is an identity in $Y_i(z^*, e^*)$ for all values of $(z^*, e^*)$, we have

$$\begin{aligned}
\forall i = 1, \ldots, n, \sum_{\mathbf{z}\in\Omega_i(z,e)} w_i(\mathbf{z}) P(\mathbf{Z}) &= \frac{1}{n} \\
\forall i = 1, \ldots, n, \sum_{\mathbf{z}\in\Omega_i(z',e')} w_i(\mathbf{z}) P(\mathbf{Z}) &= -\frac{1}{n} \\
\forall i, \forall (z^*,e^*) \neq (z,e),(z',e'), \sum_{\mathbf{z}\in\Omega_i(z^*,e^*)} w_i(\mathbf{z}) P(\mathbf{Z}) &= 0
\end{aligned}$$

Let us now show that $0 < \pi_i(z, e) < 1$ is necessary for unbiasedness. Suppose there exists a $j$ such that $\pi_j(z, e) = \sum_{\mathbf{z} \in \Omega_j(z,e)} P(\mathbf{Z}) = 0$, then $P(\mathbf{Z}) = 0 \ \forall \ \mathbf{z} \in \Omega_j(z, e)$. This means that $E[\hat{\theta}]$ is free of $Y_j(z, e)$ irrespective of $w_j(\mathbf{z})$, see line 4 of the previous equation, and hence cannot be equal to $\sum_{i=1}^n Y_i(z, e)$. Similarly, suppose there exists a $j$ such that $\pi_j(z, e) = \sum_{\mathbf{z} \in \Omega_j(z,e)} P(\mathbf{Z}) = 1$. This implies that $\Omega_j(z, e) = \Omega$. Since for fixed $j$, the sets $\Omega_j(z, e)$ are disjoint, we have $P(\mathbf{Z}) = 0$ for any $\mathbf{z} \in \Omega_j(z', e')$. Hence by the previous argument, $E[\hat{\theta}_1]$ will be free of $Y_j(z', e')$ and therefore cannot be unbiased. A similar argument will show the necessity of $0 < \pi_i(z', e') < 1$. □

### *B.2. Proof of Theorem 3.2*

*Proof.* Note that $\hat{\theta}_2$ given by 3 is contained in the class of estimators given by $\hat{\theta}$ 2, since $w(\mathbf{z}) = w(z^*, e^*)$. Using the results from Theorem 3.1, we have $\hat{\theta}_2$ is unbiased iff for each $i = 1, \ldots, n$

$$\begin{aligned}
\sum_{\mathbf{z} \in \Omega_i(z,e)} w_i(\mathbf{z}) P(\mathbf{Z}) &= \frac{1}{n} \\
\implies \sum_{\mathbf{z} \in \Omega_i(z,e)} w_i(z, e) P(\mathbf{Z}) &= \frac{1}{n} \\
\implies w_i(z, e) \sum_{\mathbf{z} \in \Omega_i(z,e)} P(\mathbf{Z}) &= \frac{1}{n}, \\
\implies w_i(z, e) \pi_i(z, e) &= \frac{1}{n} \\
\implies w_i(z, e) &= \frac{1}{n \pi_i(z, e)}
\end{aligned}$$

A similar argument shows that $w_i(z', e') = \frac{1}{n \pi_i(z', e')}$ and $w_i(z, e) = 0$ for all $(z^*, e^*) \neq (z, e)$ and $(z', e')$. □

## Appendix C: Missing proofs from Section 4

### C.1. Proof of Lemma 4.1

*Proof.* Let $0 < k \leq 1$ be a constant to be specified later and let $\hat{\theta}_1 = (1-k)\hat{\theta}$. Then we have

$$\begin{aligned} MSE(\hat{\theta}_1) &= \mathbb{E}\left((1-k)\hat{\theta} - \theta\right)^2 \\ &= \mathbb{E}\left((\hat{\theta} - \theta)^2\right) + k^2\mathbb{E}(\hat{\theta}^2) - 2k\mathbb{E}\left(\hat{\theta}^2 - \hat{\theta}\theta\right) \\ &= MSE(\hat{\theta}) + k^2(Var(\hat{\theta}) + \theta^2) - 2k(\mathbb{E}(\hat{\theta}^2) - \theta^2) \\ &= MSE(\hat{\theta}) + k^2\left(Var(\hat{\theta}) + \theta^2\right) - 2kVar(\hat{\theta}) \end{aligned} \tag{4}$$

Now if $k > 0$ and

$$k\left(Var(\hat{\theta}) + \theta^2\right) < 2Var(\hat{\theta}),$$

for all $\theta$, then from equation 4, it follows that $MSE(\hat{\theta}_1) < MSE(\hat{\theta})$. We need to show that such a $k$ exists.

To see that such a $k$ exists, first note that the MSE and variance is a function of the design $\mathbb{P}$ and the unknown but fixed potential outcomes $\{Y_i(z_i, e_i)\}_{i=1}^n$. Let us denote this set of fixed potential outcomes by $\mathbb{T}$. Next, let

$$k_0 = \min_{\mathbb{T}} \frac{2Var(\hat{\theta})}{Var(\hat{\theta}) + \theta^2}$$

Note that $k_0$ depends only on the design $\mathbb{P}$, in particular, it does not depend on the potential outcomes (some of which are unobservable) as the minimum is taken over all possible potential outcomes. Hence, $k_0$ can be computed. Since $\min_{\mathbb{T}} Var(\hat{\theta}) > 0$, we have $k_0 \geq 0$. Also, since both the numerator and denominator are strictly positive, and $\theta^2$ is finite and bounded, the denominator is not infinite, we have $k_0 > 0$. Let $k = \min(k_0, 1)$. Hence we have $0 < k \leq 1$. Now we have two cases:

Case 1: When $k_0 < 1$, $k = k_0$ and by definition of $k_0$ , $MSE(\hat{\theta}_1) < MSE(\hat{\theta})$.

Case 2: When $k_0 \geq 1$, $k = 1$, and $\hat{\theta}_1 = 0$. But $k_0 \geq 1$ implies that $2Var(\hat{\theta}) > Var(\hat{\theta}) + \theta^2$ or $Var(\hat{\theta}) \geq \theta^2$. Using this fact and substituting $k = 1$ in equation 4, one can see that $MSE(\hat{\theta}_1) < MSE(\hat{\theta})$. Note that in such a case, the variance of $\hat{\theta}$ is so large that the constant estimator 0 is able to beat it. This happens when $Var(\hat{\theta}) > \theta^2$, making estimation impossible.

In both cases, we have $MSE(\hat{\theta}_1) < MSE(\hat{\theta})$, completing the proof. ◻

### C.2. *Proof of Theorem 4.1*

*Proof of Theorem 4.1.* To show that the Horvitz-Thompson estimator is inadmissible, from Lemma 4.1, it suffices to show that the variance of the HT estimator can never be zero for a random exposure design $\mathbb{P}$. Let $U_i = I(Z_i = z', E_i = e')$ and $V_i = I(Z_i = z, E_i = e)$ and $p_i = \mathbb{E}(U_i)$ and $q_i = \mathbb{E}(V_i)$. Let us assume to the contrary that the variance of $\hat{\theta}_{HT} = 0$ for some $\theta$ for a non-constant design $\mathbb{P}$. The variance of $\hat{\theta}_{HT}$ is 0 iff

$$\begin{aligned} \hat{\theta}_{HT} &= \mathbb{E}\left(\hat{\theta}_{HT}\right) = \theta \text{ a.s. } \mathbb{P} \\ \iff \sum_{i=1}^{n}\left(Y_i(z,e)\frac{U_i}{p_i} - Y_i(z',e')\frac{V_i}{q_i}\right) &= \sum_{i=1}^{n}\left(Y_i(z,e) - Y_i(z',e')\right) \text{ a.s. } \mathbb{P} \\ \iff \sum_{i=1}^{n}\frac{Y_i(z,e)}{p_i}\left(U_i - p_i\right) &= \sum_{i=1}^{n}\frac{Y_i(z',e')}{q_i}\left(V_i - q_i\right) \text{ a.s. } \mathbb{P} \end{aligned} \tag{5}$$

Since the potential outcomes are fixed, they cannot be functions of random variables. Hence the following are the only possible solutions of equation 5:

1. For any constants $c_1, c_2$, $Y_i(z,e) = c_1 p_i$, $Y_i(z',e') = c_2 q_i$ for all $i$ and $c_1 \sum_i (U_i - p_i) = c_2 \sum_i (V_i - q_i)$ (or)
2. For any constant $c_1, Y_i(z,e) = c_1 p_i, Y_i(z',e') = 0$ for all $i$ and $c_1 \sum_i (U_i - p_i) = 0$.
3. For any constant $c_2, Y_i(z,e) = 0, Y_i(z',e') = c_2 q_i$, for all $i$ and $c_2 \sum_i (V_i - q_i) = 0$.
4. $Y_i(z,e)$ and $Y_i(z',e')$ are all 0.

Let $n_{d_0} = \sum_i U_i$ and $n_{d_1} = \sum_i V_i$. Then $\mathbb{E}[n_{d_1}] = \sum_i p_i$ and $\mathbb{E}[n_{d_0}] = \sum_i q_i$. Ignoring the last trivial solution of equation 5, we have:

1. $c_1(n_{d_1} - \mathbb{E}[n_{d_1}]) = c_2(n_{d_0} - \mathbb{E}[n_{d_0}])$ for any constants $c_1, c_2$.
2. $c_1(n_{d_1} - \mathbb{E}[n_{d_1}]) = 0$ for any constant $c_1$
3. $c_2(n_{d_0} - \mathbb{E}[n_{d_0}]) = 0$ for any constant $c_2$.

Since these hold for any $c_1$ and $c_2$, all three solutions imply that $n_{d_1} = \mathbb{E}[n_{d_1}]$ or $n_{d_0} = \mathbb{E}[n_{d_0}]$ a.s. $\mathbb{P}$. This implies that $Var(n_{d_0}) = 0$ or $Var(n_{d_1}) = 0$, which further implies that either $n_{d_0}$ or $n_{d_1}$ is constant. This contradicts the assumption that $\mathbb{P}$ is a random exposure design. □

### C.3. Proof of Theorem 4.2

The proof is a generalization of Godambe and Joshi (1965), Theorem 8.1, with several modifications to account for the potential outcomes framework. Recall that the causal effect is a contrast between potential outcomes at two different treatment and exposure combinations $d_0$ and $d_1$ defined as

$$\theta = \frac{1}{n}\left(\sum_{i=1}^{n} Y_i(d_1) - Y_i(d_0)\right) \tag{6}$$

where $d_1 = (z, e)$ and $d_0 = (z', e')$. Let $\hat{\theta} = (1/n)\left(\hat{\theta}(d_1) - \hat{\theta}(d_0)\right)$ denote the H-T estimator, where

$$\hat{\theta}(d_1) = \sum_{i=1}^{n} \frac{I(f(\mathbf{Z}) = d_1)}{\pi_i(d_1)} Y_i(d_1) \text{ and } \hat{\theta}(d_0) = \sum_{i=1}^{n} \frac{I(f(\mathbf{Z}) = d_0)}{\pi_i(d_0)} Y_i(d_0)$$

Let $\hat{\beta}$ denote any linear estimator, that is, $\hat{\beta} = (1/n)\left(\hat{\beta}(d_1) - \hat{\beta}(d_0)\right)$, where

$$\hat{\beta}(d_j) = \sum_{i=1}^{n} w_i(\mathbf{Z}, d_j) Y_i(d_j) \text{ for } j \in \{0, 1\}$$

Let $P(\mathbf{Z})$ be any fixed exposure design, i.e. $\sum_{i=1}^{n} I(f(\mathbf{Z}) = d_1) = n_{d_1}$ and $\sum_{i=1}^{n} I(f(\mathbf{Z}) = d_0) = n_{d_0}$, where $n_{d_1}$ and $n_{d_0}$ are deterministic.

Assume that the H-T estimator is inadmissible with respect to the MSE in the class of all linear estimators. Then there exists a linear estimator $\hat{\beta}$, such that

$$\mathbb{E}_{P(\mathbf{Z})}\left[(\hat{\beta} - \theta)^2\right] \leq \mathbb{E}_{P(\mathbf{Z})}\left[(\hat{\theta} - \theta)^2\right] \tag{7}$$

for all $\theta$, with strict inequality holding for at least one $\theta$. This is equivalent to equation 14 holding for all $\{Y_i(d_1), Y_i(d_0)\}_{i=1}^{n}$, with strict inequality for at least one value of $\{Y_i(d_1), Y_i(d_0)\}_{i=1}^{n}$.

Suppose

$$\hat{\beta}(d_j) = \hat{\theta}(d_j) + h_{d_j}(\mathbf{Z}), j \in \{0, 1\}, \tag{8}$$

where

$$h_{d_j}(\mathbf{Z}) = \sum_{i=1}^{n} \alpha_i(\mathbf{Z}, d_j) Y_i(d_j), \tag{9}$$

where $\alpha_i(\mathbf{Z}, d_j) = w_i(\mathbf{Z}, d_j) - b_i(d_j)$ and $b_i(d_j) = \frac{I(f(\mathbf{Z})=d_j)}{\pi_i(d_j)}$. Our goal is to show that $\alpha_i(\mathbf{Z}, d_j)$ is 0 for all $i = \{1, \ldots, n\}$, for all $\mathbf{Z}$ such that $P(\mathbf{Z}) > 0$, and $j = \{0, 1\}$.

Substituting the values of $\hat{\beta}(d_j)$ from equation 15 in equation 14, we get:

$$\mathbb{E}_{P(\mathbf{Z})}\left[\left(\left(\hat{\theta}(d_1) + h(d_1, \mathbf{Z}) - \theta(d_1)\right) - \left(\hat{\theta}(d_0) + h(d_0, , \mathbf{Z}) - \theta(d_0)\right)\right)^2\right] \leq \mathbb{E}_{P(\mathbf{Z})}\left[\left(\left(\hat{\theta}(d_1) - \theta(d_1)\right) - \left(\hat{\theta}(d_0) - \theta(d_0)\right)\right)^2\right] \quad (10)$$

Algebraic manipulation gives us

$$\mathbb{E}_{P(\mathbf{Z})}\left[h_{d_1}(\mathbf{Z})^2 + 2h_{d_1}(\mathbf{Z}) \cdot \left(\hat{\theta}(d_1) - \theta(d_1)\right) + h_{d_0}(\mathbf{Z})^2 + 2h_{d_0}(\mathbf{Z}) \cdot \left(\hat{\theta}(d_0) - \theta(d_0)\right)\right] \leq \mathbb{E}_{P(\mathbf{Z})}\left[2\left(\hat{\theta}(d_0) - \theta(d_0)\right) h_{d_1}(\mathbf{Z}) + 2\left(\hat{\theta}(d_1) - \theta(d_1)\right) h_{d_0}(\mathbf{Z})\right] \quad (11)$$

Substituting the definition of $h_{d_j}(\mathbf{Z})$ from equation 16 into the above inequality, and carefully inspecting, one can see that equation 18 is a semi-negative definite quadratic form in $Y_i(d_1)$ and $Y_i(d_0)$. Hence the coefficients of $Y_i(d_1)^2$ and $Y_i(d_0)^2$ are negative for each $i$. Note that these terms appear only on the LHS of inequality 18, as the RHS involves products of form $Y_i(d_0)Y_i(d_1)$, and hence can be ignored. Extracting the coefficients of $Y_i(d_1)^2$ and $Y_i(d_0)^2$ and setting them to be less than 0, we get

$$\mathbb{E}_{P(\mathbf{Z})}\left[\alpha_i(\mathbf{Z}, d_j)^2 + 2\alpha_i(\mathbf{Z}, d_j)(b_i(d_j) - 1)\right] \leq 0, \quad (12)$$

for all $i = \{1, \ldots, n\}$, $j = \{0, 1\}$, and $\mathbf{Z} \in \Omega_i(d_j)$, where $\Omega_i(d_j)$ is the set of treatment assignments that reveals potential outcome $Y_i(d_j)$ for unit $i$. Now, let $\mathbb{E}_{P(\mathbf{Z})}[\alpha_i(\mathbf{Z}, d_j)] = \delta_i(d_j)$, then

$$\begin{aligned}
&\mathbb{E}_{P(\mathbf{Z})}\left[\left(\alpha_i(\mathbf{Z}, d_j) - \delta_i(d_j)\right)^2\right] \\
&= \mathbb{E}_{P(\mathbf{Z})}\left[\alpha_i(\mathbf{Z}, d_j)^2\right] + \mathbb{E}_{P(\mathbf{Z})}\left[\delta_i(d_j)^2\right] - 2\mathbb{E}_{P(\mathbf{Z})}\left[\alpha_i(\mathbf{Z}, d_j) \cdot \delta_i(d_j)\right] \\
&= \mathbb{E}_{P(\mathbf{Z})}\left[\alpha_i(\mathbf{Z}, d_j)^2\right] - \mathbb{E}_{P(\mathbf{Z})}\left[\delta_i(d_j)^2\right] \\
&= \mathbb{E}_{P(\mathbf{Z})}\left[\alpha_i(\mathbf{Z}, d_j)^2\right] - \delta_i(d_j)^2 \mathbb{E}_{P(\mathbf{Z})}\left[\mathbb{I}(\mathbf{Z} \in \Omega_i(d_j))\right] \\
&= \mathbb{E}_{P(\mathbf{Z})}\left[\alpha_i(\mathbf{Z}, d_j)^2\right] - \delta_i(d_j)^2 \pi_i(d_j)
\end{aligned}$$

where the last equation follows from the definition of $\pi_i(d_j)$. Hence,

$$\mathbb{E}_{P(\mathbf{Z})}\left[\alpha_i(\mathbf{Z},d_j)^2\right] = \mathbb{E}_{P(\mathbf{Z})}\left[\left(\alpha_i(\mathbf{Z},d_j)-\delta_i(d_j)\right)^2\right] + \delta_i(d_j)^2\cdot\pi_i(d_j).$$

Substituting this in equation 20, we get

$$\begin{aligned}
&\mathbb{E}_{P(\mathbf{Z})}\left[\left(\alpha_i(\mathbf{Z},d_j)-\delta_i(d_j)\right)^2\right] + \delta_i(d_j)^2\pi_i(d_j) + 2\mathbb{E}\left[\alpha_i(\mathbf{Z},d_j)(b_i(d_j)-1)\right] \le 0\\
&\implies \mathbb{E}_{P(\mathbf{Z})}\left[\left(\alpha_i(\mathbf{Z},d_j)-\delta_i(d_j)\right)^2\right] + \delta_i(d_j)^2\pi_i(d_j) + 2(b_i(d_j)-1)\mathbb{E}_{P(\mathbf{Z})}\left[\alpha_i(\mathbf{Z},d_j)\right] \le 0\\
&\implies \mathbb{E}_{P(\mathbf{Z})}\left[\left(\alpha_i(\mathbf{Z},d_j)-\delta_i(d_j)\right)^2\right] + \delta_i(d_j)^2\pi_i(d_j) + 2(b_i(d_j)-1)\delta_i(d_j) \le 0
\end{aligned} \tag{13}$$

for all $i = 1,\dots,n$, $j = \{0,1\}$ and $\mathbf{Z} \in \Omega_i(d_j)$. Note that for $\mathbf{Z} \in \Omega_i(d_j)$, $b_i(d_j) = \frac{1}{\pi_i} \ge 1$. Hence, inequality 13 implies that $\delta_i(d_j) \le 0$.

Now, consider the causal effect $\theta$ at the following specific values of the potential outcomes: $Y_i^*(d_1) = \pi_i(d_1)$, and $Y_i^*(d_0) = 0$. Hence $\theta = \sum_i \pi_i(d_1) = n_{d_1}$, by the assumption of the fixed exposure design. Also, $\hat{\theta} = \sum_{i=1}^n I(f(\mathbf{Z} = d_1)) = n_{d_1}$. Thus for $\{Y_i^*(d_1), Y_i^*(d_0)\}$, we have $\theta = \hat{\theta} = n_{d_1}$. From equation 14, this implies that $\hat{\beta} = n_{d_1} = \hat{\beta}(d_1)$. Hence, since $\hat{\beta}(d_1) = \hat{\theta}(d_1)$, we have, $h_{d_1}(\mathbf{Z}) = 0$ at $\{Y_i^*(d_1), Y_i^*(d_0)\}$, which, by definition of $h_{d_1}$ implies that $\sum_{i=1}^n \alpha_i(\mathbf{Z},d_1)\pi_j(d_j) = 0$. This implies that $\alpha_i(\mathbf{Z},d_1) = 0$ for all $i$ such that $\pi_i(d_1) > 0$. Thus, we have $w_i(\mathbf{Z},d_1) = b_i(d_1)$. A similar argument using $Y_i^*(d_1) = 0$ and $Y_i^*(d_0) = \pi_i(d_0)$ shows that $w_i(\mathbf{Z},d_0) = b_i(d_0)$ for all $i$ such that $\pi_i(d_0) > 0$, proving that $\hat{\beta} = \hat{\theta}$.

## Appendix D: Missing proofs from Section 5

### *D.1. Proof of Lemma 5.1*

*Proof.* Recall that $\theta = \frac{1}{n}\sum_{i=1}^n (Y_i(d_1) - Y_i(d_0))$. where $d_1 = (z,e)$ and $d_0 = (z',e')$. Let $\hat{\theta} = \hat{\theta}(d_1) - \hat{\theta}(d_0)$ , where

$$\hat{\theta}(d_1) = (1/n)\cdot\sum_{i=1}^n w_i(\mathbf{Z},d_1)Y_i(d_1) \text{ and } \hat{\theta}(d_0) = (1/n)\cdot\sum_{i=1}^n w_i(\mathbf{Z},d_0)I(f(\mathbf{Z}) = d_0)Y_i(d_0)$$

Let $\hat{\theta}'$ denote any linear estimator, that is, $\hat{\theta}' = \hat{\theta}'(d_1) - \hat{\theta}'(d_0)$, where

$$\hat{\theta}'(d_j) = (1/n)\cdot\sum_{i=1}^n w_i'(\mathbf{Z},d_j)Y_i(d_j) \text{ for } j \in \{0,1\}$$

Note that $w_i(\mathbf{Z}, d_j)$ and $w'_i(\mathbf{Z}, d_j)$ are nonzero iff $I(D_j = d_j) = 1$, i.e. the estimators are functions of observed potential outcomes. Assume that $\hat{\theta}$ is inadmissible with respect to the MSE in the class of all linear estimators. Then there exists an estimator $\hat{\theta}'$ that dominates $\hat{\theta}$, i.e.

$$\mathbb{E}_{P(\mathbf{Z})}\left[(\hat{\theta}' - \theta)^2\right] \leq \mathbb{E}_{P(\mathbf{Z})}\left[(\hat{\theta} - \theta)^2\right] \tag{14}$$

for all $\theta$, with strict inequality holding for at least one $\theta$. This is equivalent to equation 14 holding for all $\{Y_i(d_1), Y_i(d_0)\}_{i=1}^n$, with strict inequality for at least one value of $\{Y_i(d_1), Y_i(d_0)\}_{i=1}^n$. Now, suppose

$$\hat{\theta}'(d_j) = \hat{\theta}(d_j) + \Delta_{d_j}(\mathbf{Z}),\ j \in \{0, 1\}, \tag{15}$$

where

$$\Delta_{d_j}(\mathbf{Z}) = \sum_{i=1}^{n} \alpha_i(\mathbf{Z}, d_j) Y_i(d_j), \tag{16}$$

where $\alpha_i(\mathbf{Z}, d_j) = w_i(\mathbf{Z}, d_j) - w'_i(\mathbf{Z}, d_j)$. Our goal is to show that $\alpha_i(\mathbf{Z}, d_j)$ is 0 for all $i = \{1, \dots, n\}$, for all $\mathbf{Z}$ such that $P(\mathbf{Z}) > 0$, and $j = \{0, 1\}$, which implies that $\hat{\theta}' = \hat{\theta}$.

Substituting the values of $\hat{\theta}'(d_j)$ from equation 15 in equation 14, we get:

$$\begin{aligned} &\mathbb{E}_{P(\mathbf{Z})}\left[\left(\left(\hat{\theta}(d_1) + \Delta_{d_1}(\mathbf{Z}) - \theta(d_1)\right) - \left(\hat{\theta}(d_0) + \Delta_{d_0}(\mathbf{Z}) - \theta(d_0)\right)\right)^2\right] \\ &\qquad \leq \mathbb{E}_{P(\mathbf{Z})}\left[\left(\left(\hat{\theta}(d_1) - \theta(d_1)\right) - \left(\hat{\theta}(d_0) - \theta(d_0)\right)\right)^2\right] \end{aligned} \tag{17}$$

Algebraic manipulation gives us

$$\begin{aligned} &\mathbb{E}_{P(\mathbf{Z})}\left[\Delta_{d_1}(\mathbf{Z})^2 + 2\Delta_{d_1}(\mathbf{Z}) \cdot \left(\hat{\theta}(d_1) - \theta(d_1)\right) + \Delta_{d_0}(\mathbf{Z})^2 + 2\Delta_{d_0}(\mathbf{Z}) \cdot \left(\hat{\theta}(d_0) - \theta(d_0)\right)\right] \\ &\qquad \leq \mathbb{E}_{P(\mathbf{Z})}\left[2\left(\hat{\theta}(d_0) - \theta(d_0)\right)\Delta_{d_1}(\mathbf{Z}) + 2\left(\hat{\theta}(d_1) - \theta(d_1)\right)\Delta_{d_0}(\mathbf{Z})\right] \end{aligned} \tag{18}$$

Substituting the definition of $\Delta_{d_j}(\mathbf{Z})$ from equation 16 into the above inequality, and carefully inspecting, one can see that equation 18 is a semi-negative definite quadratic form in $Y_i(d_1)$ and $Y_i(d_0)$. Hence the coefficients of $Y_i(d_1)^2$ and $Y_i(d_0)^2$ are negative for each $i$. Note that these terms appear only on the LHS of inequality 18, as the RHS involves products of form $Y_i(d_0)Y_i(d_1)$, and hence can be ignored.

Extracting the coefficients of $Y_i(d_1)^2$ and $Y_i(d_0)^2$ and setting them to be less than 0, we get

$$\mathbb{E}_{P(\mathbf{Z})}\left[\alpha_i(\mathbf{Z}, d_j)^2 + 2\alpha_i(\mathbf{Z}, d_j)(w_i(\mathbf{Z}, d_j) - 1)\right] \leq 0, \tag{19}$$

for all $i = \{1, \ldots, n\}$, $j = \{0, 1\}$. This implies

$$\mathbb{E}_{\mathbf{Z}\in\Omega_i(d_j)}\left[\alpha_i(\mathbf{Z}, d_j)(w_i(\mathbf{Z}, d_j) - 1)\right] \leq 0, \tag{20}$$

where $\Omega_i(d_j) = \{\mathbf{z} \in \Omega : I_i(d_j) > 0\}$. Note that by positivity assumption $\pi_i(d_j) > 0$, for each $i$, $\Omega_i(d_j)$ is non-empty, hence equation 20 is well defined for each $i$. Now let $c_i = \min_{\mathbf{z}\in\Omega_i}(w_i(\mathbf{z}, d_j) - 1)$. By assumption, $c_i > 0$. hence we have

$$0 \geq \mathbb{E}_{P(\mathbf{Z})}\left[\alpha_i(\mathbf{Z}, d_j)(w_i(\mathbf{Z}, d_j) - 1)\right] \geq c_i\mathbb{E}_{P(\mathbf{Z})}\left[\alpha_i(\mathbf{Z}, d_j)\right]$$

which implies $\mathbb{E}_{P(\mathbf{Z})}\left[\alpha_i(\mathbf{Z}, d_j)\right] \leq 0$.

Now we will show that $\sum_i \mathbb{E}[\alpha_i(\mathbf{Z}, d_j)] = 0$. Recall that for a fixed $j$, there exists a set of non zero potential outcomes $\{Y_i^*(d_j)\}$ such that $MSE(\hat{\theta}(d_j), \theta^*(d_j)) = 0$. This implies that $MSE(\hat{\theta}'(d_j), \theta^*(d_j)) = 0$, as $\hat{\theta}'$ dominates $\hat{\theta}$. In particular, it implies that at $\{Y_i^*(d_j)\}$, $Var(\hat{\theta}') = Var(\hat{\theta}) = 0$, and hence $\hat{\theta}'(d_j) = \hat{\theta}(d_j) = \theta^*(d_j)$ at $\{Y_i^*(d_j)\}$. This implies

$$0 = \hat{\theta}' - \hat{\theta} = \sum_{i=1}^{n} \alpha_i(\mathbf{Z}, d_j)Y_i^*(d_j) = 0$$

Since $Y_i^*(d_j) \neq 0$, we have $\sum_{i=1}^{n} \alpha_i(\mathbf{Z}, d_j) = 0$, and hence $\sum_{i=1}^{n} \mathbb{E}[\alpha_i(\mathbf{Z}, d_j)] = 0$. We have already shown that $\mathbb{E}[\alpha_i(\mathbf{Z}, d_j)] \leq 0$. Thus, this implies that $\mathbb{E}[\alpha_i(\mathbf{Z}, d_j)] = 0$.

Now let $\delta_i(d_j) = \mathbb{E}[\alpha_i(\mathbf{Z}, d_j)]$, and note that $\mathbb{E}[(\alpha_i(\mathbf{Z}, d_j) - \delta_i(d_j))^2] = \mathbb{E}[\alpha_i(\mathbf{Z}, d_j)^2] - \delta_i(d_j)^2$, which implies

$$\mathbb{E}[\alpha_i(\mathbf{Z}, d_j)^2] = E[(\alpha_i(\mathbf{Z}, d_j) - \delta_i(d_j))^2] + \delta_i(d_j)^2.$$

Substituting this in equation 20, we get

$$\mathbb{E}[(\alpha_i(\mathbf{Z}, d_j) - \delta_i(d_j))^2] + \delta_i(d_j)^2 + 2\mathbb{E}[\alpha_i(\mathbf{Z}, d_j)(w_i(\mathbf{Z}, d_j) - 1)] \leq 0,$$

Let $c_i = \min_z w_i(\mathbf{z}, d_j) - 1$, we get

$$\mathbb{E}[(\alpha_i(\mathbf{Z}, d_j) - \delta_i(d_j))^2] + \delta_i(d_j)^2 + 2c_i\delta_i(d_j) \leq 0,$$

But we have shown that $\delta_i(d_j) = 0$, hence we get $\mathbb{E}[\alpha_i(\mathbf{Z}, d_j)^2] \leq 0$. This implies $\alpha_i(\mathbf{Z}, d_j) = 0$, whenever $\pi_i > 0$. Hence $w_i'(\mathbf{Z}, d_j) = w_i(\mathbf{Z}, d_j)$ which implies $\hat{\theta}' = \hat{\theta}$.

□

### D.2. *Proof of Theorem 5.1*

*Proof.* Consider the set of potential outcomes $Y_i^*(d_j) = c$, then by definition,

$$\hat{\theta}_J(d_j) = \frac{1}{\hat{n}} \sum_{i=1}^{n} \frac{I(D_i = d_j)}{\pi_i(d_j)} c = c$$

and $\theta(d_j) = c$. Hence $MSE(\hat{\theta}_J(d_j), \theta(d_j)) = 0$. Moreover,

$$\hat{\theta}_J(d_j) = \frac{1}{n} \sum_{i=1}^{n} w_i(\mathbf{Z}, d_j) Y_i(d_j),$$

where $w_i(\mathbf{Z}, d_j) = \frac{n}{\hat{n}(d_j)} \frac{I(D_i=d_j)}{\pi_i(d_j)}$. We need $\min_{\mathbf{z}\in\Omega_i(d_j)} w_i(\mathbf{z}, d_j) > 1$. Since and $0 < \pi_i(d_j) < 1$, and $\Omega_i(d_j) = \{\mathbf{z} : I(d_j) > 1\}$, $\hat{n}(d_j) < n$ implies $w_i(\mathbf{z}, d_j) > 1$ for all $\mathbf{z}$ in $\Omega_i$. Thus the conditions of Lemma 5.1 are satisfied. Hence $\hat{\theta}_J$ is admissible. □

### D.3. *Proof of Theorem 5.3*

*Proof.* Consider the set of potential outcomes $Y_i^*(d_j) = \pi_i(d_j)$, then by definition,

$$\hat{\theta}_R(d_j) = \frac{1}{n} \frac{\sum_i \pi_i(d_j)}{n_{d_j}} \sum_{i=1}^{n} \frac{I(D_i = d_j)}{\pi_i(d_j)} \pi_i(d_j) = \frac{\sum_i \pi_i(d_j)}{n}$$

and $\theta(d_j) = \sum_i \pi_i(d_j)$. Hence $MSE(\hat{\theta}_R(d_j), \theta(d_j)) = 0$. From the definition, $\hat{\theta}_R(d_j) = (1/n) \sum_{i=1}^{n} w_i(\mathbf{Z}, d_j) Y_i(d_j)$, where $w_i(\mathbf{Z}, d_j) = \frac{\sum_i \pi_i(d_j)}{n_{d_j}} \frac{I(D_i=d_j)}{\pi_i(d_j)}$. Hence if $n_{d_j} < \sum_i \pi_i(d_j)$, and $0 < \pi_i(d_j) < 1$, then $w_i(\mathbf{z}, d_j) > 1$ for $\mathbf{z} \in \Omega_i$. Thus the conditions of Lemma 5.1 are satisfied. Hence $\hat{\theta}_R$ is admissible. □

### D.4. *Proof of Theorem 5.2*

Let $I_1, \ldots, I_n$ be $\{0, 1\}$–valued random variables with $\Pr(I_i = 1) = \pi_i$. Assume the dependence structure is given by an (undirected) graph $G$ on $\{1, \ldots, n\}$ and that two indices $i \neq j$ are dependent only if their graph distance is 1 or 2 (i.e. they are adjacent or share a common neighbour). Let $d_i$ denote the degree of node $i$ in $G$. Define

$$S_n = \sum_{i=1}^{n} \frac{I_i}{\pi_i}, \qquad \sigma^2 = \mathrm{Var}(S_n).$$

We begin by deriving a normal approximation to the dependent sum $S_n$ by using the classic result of Chen and Shao (2004).

**Lemma D.1.** *Suppose the following hold:*

1. *(bounded degree)* $d_{\max} := \max_{1\le i\le n} d_i < \infty$ *is a constant (independent of* $n$*);*
2. *(uniform lower bound on marginals)* $\pi_{\min} := \min_{1\le i\le n} \pi_i > 0$*;*
3. *(nondegeneracy of variance) There exists a* $c$ *(independent of* $n$*) such that* $\sigma^2 \ge cn$

*Then, writing* $\Phi$ *for the standard normal distribution function,*

$$\sup_{z\in\mathbb{R}} \left| \Pr\left( \frac{S_n - \mathbb{E}S_n}{\sigma} \le z \right) - \Phi(z) \right| \;=\; O\big(n^{-1/2}\big).$$

*Proof.* Set

$$Z_i \;:=\; \frac{I_i}{\pi_i} - 1, \qquad W := \sum_{i=1}^{n} Z_i = S_n - \mathbb{E}S_n,$$

and let $\sigma^2 = \mathrm{Var}(W)$. We will apply Chen and Shao (2004), Theorem 2.2 (display (2.6) there) to the normalized variables

$$X_i := \frac{Z_i}{\sigma},$$

so that $\mathrm{Var}\big(\sum_i X_i\big) = 1$. We begin by verifying the condition (LD2) in the Chen and Shao (2004). For each $i$ let $A_i$ be the closed two–hop neighbourhood of $i$ in $G$,

$$A_i := \{j : \mathrm{dist}_G(i, j) \le 2\}.$$

Take $B_i := A_i$. By construction $X_i$ is independent of $\{X_j : j \notin A_i\}$, and $X_{A_i}$ is independent of the complement of $B_i$; hence the local dependence condition (LD2) in Chen and Shao (2004) holds with these sets.

Next, we bound the parameters $\kappa$ and $\theta$. Define $Y_i := \sum_{j\in A_i} X_j$ (Chen–Shao notation). Because $I_j/\pi_j - 1$ is bounded in absolute value by $1/\pi_j \le 1/\pi_{\min}$, we have for every $j$

$$|X_j| = \frac{|I_j/\pi_j - 1|}{\sigma} \le \frac{1}{\pi_{\min}\,\sigma}.$$

Under bounded degree $d_{\max} = O(1)$ the two-hop neighborhoods are uniformly bounded: there is a constant $C_A$ (depending only on $d_{\max}$) with

$$|A_i| \le C_A \quad \text{for all } i.$$

Therefore

$$|Y_i| \leq |A_i| \cdot \frac{1}{\pi_{\min}\,\sigma} \leq \frac{C_A}{\pi_{\min}\,\sigma}.$$

Set

$$\theta^p \;:=\; \max_i \left(\mathbb{E}|X_i|^p + \mathbb{E}|Y_i|^p\right).$$

For any fixed $p \geq 2$ (and in particular for $p = 4$) the boundedness above gives

$$\theta \;\leq\; \frac{C}{\pi_{\min}\,\sigma}$$

for an absolute constant $C$ (depending on $C_A$ and $p$). Since by assumption (c) and the variance lower bound proved earlier we have $\sigma^2 \geq cn$ for some $c > 0$, it follows that

$$\theta = O\big(\sigma^{-1}\big) = O\big(n^{-1/2}\big).$$

Define the overlap parameter

$$\kappa \;:=\; \max_i \big|\{\, j : A_j \cap A_i \neq \varnothing \,\}\big|.$$

Because $|A_i| \leq C_A = O(1)$ and each $A_j$ is itself of size $O(1)$, the number of indices $j$ whose neighborhood intersects $A_i$ is $O(1)$. Hence $\kappa = O(1)$. Chen–Shao give an explicit Kolmogorov bound (their display (2.6)); specializing to $p = 4$ and inserting the above bounds yields (for absolute constants $C_1, C_2$)

$$\sup_z \big|F(z) - \Phi(z)\big| \;\leq\; C_1\,\kappa\,n\,\theta^3 \;+\; C_2\,\theta^2\sqrt{\kappa n},$$

where $F$ is the cdf of $\sum_i X_i$. Since $\kappa = O(1)$ and $\theta = O(n^{-1/2})$, both $n\theta^3$ and $\theta^2\sqrt{n}$ are $O(n^{-1/2})$. Thus the right-hand side is $O(n^{-1/2})$, and therefore

$$\sup_z \left|\Pr\left(\frac{W}{\sigma} \leq z\right) - \Phi(z)\right| = O(n^{-1/2}).$$

Returning to the original (un-normalized) sum $S_n$, we have shown

$$\sup_z \left|\Pr\left(\frac{S_n - \mathbb{E}S_n}{\sigma} \leq z\right) - \Phi(z)\right| = O\big(n^{-1/2}\big),$$

as claimed. □

Next, we consider an set of indicators $\{I_i^*\}$ where each $I_i^*$ is independent, and has the same marginal distribution as $I_i$. Define the sum $S_n^* = \sum_i I_i^*/\pi_i$ and note that $E[S_n] = E[S_n^*] = n$. We derive an Edgeworth expansion for $S_n^*$. The following lemma gives the moments of $S_n^*$ that are needed for the Edgeworth expansion

**Lemma D.2** (Moments of the HT count $\widehat{n}^*$)**.** *Let $I_1^*, \dots, I_n^*$ be independent Bernoulli random variables with $P(I_i^* = 1) = \pi_i \in (0, 1)$. Define*

$$X_i^* := \frac{I_i^*}{\pi_i}, \qquad \widehat{n}^* := \sum_{i=1}^{n} X_i^*.$$

*Then for each i*

$$E[X_i^*] = 1, \qquad E[(X_i^*)^k] = \pi_i^{1-k} \;\; (k \geq 1),$$

*and therefore*

$$E[\widehat{n}^*] = n,$$

$$\mathrm{Var}(\widehat{n}^*) = \sum_{i=1}^{n} \frac{1-\pi_i}{\pi_i},$$

$$\mu_3(\widehat{n}^*) := E\left[(\widehat{n}^* - n)^3\right] = \sum_{i=1}^{n} \frac{(1-\pi_i)(1-2\pi_i)}{\pi_i^2}.$$

*Consequently the standardized (aggregate) skewness is*

$$\gamma \;=\; \frac{\mu_3(\widehat{n}^*)}{\left(\mathrm{Var}(\widehat{n}^*)\right)^{3/2}}.$$

*Proof.* All steps are elementary algebraic calculations. First, we compute the raw moments of $X_i^*$. Since $X_i^*$ takes the value $r_i := 1/\pi_i$ with probability $\pi_i$ and 0 with probability $1 - \pi_i$,

$$E[(X_i^*)^k] = r_i^k \cdot \pi_i + 0 \cdot (1-\pi_i) = \frac{1}{\pi_i^k}\, \pi_i = \pi_i^{1-k}.$$

In particular $E[X_i^*] = \pi_i^0 = 1$. Linearity of expectation gives

$$E[\widehat{n}^*] = \sum_{i=1}^{n} E[X_i^*] = \sum_{i=1}^{n} 1 = n.$$

Now we compute the variance of $\widehat{n}^*$. First compute $\mathrm{Var}(X_i^*) = E[(X_i^*)^2] - E[X_i^*]^2$. Using $E[(X_i^*)^2] = \pi_i^{-1}$ and $E[X_i^*] = 1$,

$$\mathrm{Var}(X_i^*) = \frac{1}{\pi_i} - 1 = \frac{1-\pi_i}{\pi_i}.$$

Hence,

$$\mathrm{Var}(\widehat{n}) = \sum_{i=1}^{n} \mathrm{Var}(X_i^*) = \sum_{i=1}^{n} \frac{1-\pi_i}{\pi_i}.$$

The third central moment is computed as follows. For each $i$,

$$E[(X_i^* - 1)^3] = E[(X_i^*)^3] - 3E[(X_i^*)^2] + 3E[X_i^*] - 1.$$

Using $E[(X_i^*)^3] = \pi_i^{-2},\ E[(X_i^*)^2] = \pi_i^{-1},\ E[X_i^*] = 1$, we obtain

$$E[(X_i^* - 1)^3] = \frac{1}{\pi_i^2} - \frac{3}{\pi_i} + 2 = \frac{(1-\pi_i)(1-2\pi_i)}{\pi_i^2}.$$

Independence implies the third central moment of the sum is the sum of individual third central moments:

$$\mu_3(\widehat{n}^*) = \sum_{i=1}^{n} E[(X_i^* - 1)^3] = \sum_{i=1}^{n} \frac{(1-\pi_i)(1-2\pi_i)}{\pi_i^2}.$$

By definition, the standardized skewness $\gamma = \mu_3(\widehat{n}^*)/(\mathrm{Var}(\widehat{n}^*))^{3/2}$, giving the stated formula. This completes the proof. □ □

We now derive the Edgeworth expansion of $S_n^*$ using the results of Angst and Poly (2017). Since $S_n^*$ is a sum of discrete random variables, we need a Edgeworth expansion for non-lattice discrete random variables. Angst and Poly (2017) show that if the summands satisfy the so called *mean weak cramer condition*, then one can derive an Edgeworth expansion for discrete non-lattice random variables. The next lemma shows that the summands of $S_n^*$ indeed satisfy the mean weak Cramer condition.

**Lemma D.3** (Block-sum reduction and mean weak Cramér property)**.** *Let $I_i^*$ be independent* Bernoulli$(\pi_i)$ *random variables and define*

$$X_i^* = \frac{I_i^*}{\pi_i}, \qquad i = 1, \ldots, n.$$

*Assume that the sequence* $\{\pi_i\}$ *is a sequence of* generic *real numbers between* $(0,1)$ *and there exists a constant* $a > 0$ *such that* $\pi_i > a$. *Let* $S_n^* = \sum_i X_i^*$. *Then there exist* $Y_1, \ldots, Y_K$ *such that* $S_n^* = \sum_k Y_k + o(n)$ *and the sequence* $Y_m$ *satisfies the mean weak Cramer condition: there exist constants* $C > 0$, $b > 0$, *and* $R > 0$ *such that for all* $|t| > R$ *and all sufficiently large* $m$*,*

$$\frac{1}{m}\sum_{k=1}^{m} |\varphi_{Y_k}(t)| \leq 1 - \frac{C}{|t|^b}.$$

*Consequently,* $S_n^*$ *admits a one-term Edgeworth expansion given by for the standardized sum* $Z_n^* = S_n^*/\sigma$*,*

$$P(Z_n^* \leq x) = \Phi(x) + \frac{\gamma}{6}(1 - x^2)\phi(x) + O(1/n),$$

*Proof.* The construction follows example 3.1 in Angst and Poly (2017). Partition the indices $\{1, \ldots, n\}$ into $m = \lfloor n/3 \rfloor$ disjoint blocks $B_1, \ldots, B_m$ of size 3. Let

$$Y_k := \sum_{i \in B_k} X_i^*, \qquad k = 1, \ldots, m,$$

Let $L$ be the (at most two) leftover indices. Then

$$S_n^* = \sum_{k=1}^{m} Y_k + \sum_{i \in L} X_i^* = \sum_{k=1}^{m} Y_k + o(n),$$

since, for each leftover index $i$, $0 \leq X_i^* \leq 1/\pi_i \leq 1/a$, hence $|R_n| := \left|\sum_{i \in L} X_i^*\right| \leq 2/a = O(1) = o(n)$. Note that each block law has at least three atoms: Fix a block $B_k = \{i_1, i_2, i_3\}$. Each $X_{i_j}^*$ takes values in $\{0, 1/\pi_{i_j}\}$, so the block sum $Y_k = X_{i_1}^* + X_{i_2}^* + X_{i_3}^*$ takes values among the $2^3 = 8$ possible sums

$$\sum_{j \in J} \frac{1}{\pi_{i_j}}, \qquad J \subseteq \{1, 2, 3\}.$$

Having fewer than three distinct values would require special equalities among the $\{1/\pi_{i_j}\}$, such as $1/\pi_{i_1} = 1/\pi_{i_2}$ or $1/\pi_{i_1} + 1/\pi_{i_2} = 1/\pi_{i_3}$. Under the genericity assumption, such coincidences do not occur. Hence, each $Y_k$ has at least three distinct atoms with positive probabilities.

By Proposition 2.4 of Angst and Poly (2017), any family of discrete independent random variables with at least three generic atoms obeys the *mean weak Cramér* *condition*. Since the $\{Y_k\}$ are independent (their blocks are disjoint) and each has at least three generic atoms, the conditions of the proposition are satisfied. Thus there exist constants $C > 0$, $b > 0$, and $R > 0$ such that

$$\frac{1}{m}\sum_{k=1}^{m} |\varphi_{Y_k}(t)| \leq 1 - \frac{C}{|t|^b}, \qquad |t| > R.$$

Combining $S_n^* = \sum_{k=1}^{m} Y_k + R_n$ with $R_n = o(n)$ completes the first part of the proof. Next, by Theorem 4.3 of Angst and Poly (2017), the sum $S_n^*$ admits an Edgeworth expansion if (a) $Y_k$ satisfy the mean weak cramer condition (b) the smallest eigen value of $1/n \sum_k \mathrm{Var}(Y_k)$ is bounded away from 0 uniformly in $n$, and (c) Let $\rho_s(n) = (1/n)\sum_k E[|Y_k|^s]$, then there exists an $s \geq 3$ such that $\rho_s(n)$ is bounded and $\lim_{n\to\infty} (1/n)\sum_j E[\mathbf{1}\{|Y_k| > \epsilon\sqrt{n}\}|Y_k|^s] = 0$. We showed (a), and (b) and (c) automatically hold for $Y_k$ for $s = 4$, as they are bounded discrete random variables. Thus, $S_n^*$ admits a classic Edgeworth expansion with error term $O(1/n)$. □

**Lemma D.4** (Transfer of the first–order Edgeworth term)**.** *Let* $I_1, \ldots, I_n$ *be* $\{0, 1\}$*–valued random variables with* $\Pr(I_i = 1) = \pi_i$*, and suppose their dependence structure is determined by a graph* $G$ *on* $\{1, \ldots, n\}$ *such that two indices* $i \neq j$ *are dependent only if their graph distance is at most* 2*. Define*

$$S_n = \sum_{i=1}^{n} \frac{I_i}{\pi_i}, \qquad S_n^* = \sum_{i=1}^{n} \frac{I_i^*}{\pi_i},$$

*where the* $I_i^*$ *are independent Bernoulli*$(\pi_i)$ *variables. Write* $\mu_n = \mathbb{E}S_n = \mathbb{E}S_n^*$ *and* $\sigma_n^2 = \mathrm{Var}(S_n)$*. Assume:*

1. *the degrees of* $G$ *are uniformly bounded,* $d_{\max} = O(1)$*;*
2. *the probabilities* $\pi_i$ *are uniformly bounded away from* 0*,* $\pi_{\min} > 0$*;*
3. *the variance is nondegenerate,* $\sigma_n^2 \geq cn$ *for some constant* $c > 0$*.*

*Suppose moreover that the independent sum* $S_n^*$ *satisfies a first–order Edgeworth expansion at its mean,*

$$\Pr\big(S_n^* < \mu_n\big) \; = \; \tfrac{1}{2} \; + \; C\,\gamma_n^* \; + \; r_n, \tag{21}$$

*where*

$$\gamma_n^* \; = \; \frac{\kappa_3(S_n^*)}{(\sigma_n^*)^3}$$

*is the standardized third cumulant (skewness) of* $S_n^*$*, and the remainder* $r_n = o(n^{-1/2})$*. Then under (a)–(c),*

$$\Pr(S_n < \mu_n) \; = \; \tfrac{1}{2} \; + \; C\,\gamma_n^* \; + \; O(n^{-1/2}). \tag{22}$$

*Proof.* From Lemma D.1, the normalized dependent sum satisfies a Berry–Esseen bound of order $n^{-1/2}$:

$$\sup_{z\in\mathbb{R}}\left|\Pr\left(\frac{S_n-\mu_n}{\sigma_n}\leq z\right)-\Phi(z)\right|=O(n^{-1/2}).$$

For the independent sum $S_n^*$, the classical Berry–Esseen theorem gives the same rate. Hence

$$\begin{aligned}\Delta_n &:= \left|\Pr(S_n<\mu_n)-\Pr(S_n^*<\mu_n)\right|\\ &\leq \sup_z\left|\Pr\left(\frac{S_n-\mu_n}{\sigma_n}\leq z\right)-\Phi(z)\right|+\sup_z\left|\Pr\left(\frac{S_n^*-\mu_n}{\sigma_n^*}\leq z\right)-\Phi(z)\right|\\ &= O(n^{-1/2}).\end{aligned}$$

Using the Edgeworth expansion (21) for $S_n^*$,

$$\Pr(S_n^*<\mu_n)=\tfrac{1}{2}+C\,\gamma_n^*+r_n,\qquad r_n=o(n^{-1/2}),$$

and adding the bound on $\Delta_n$ gives

$$\begin{aligned}\Pr(S_n<\mu_n)&=\Pr(S_n^*<\mu_n)+O(n^{-1/2})\\ &=\tfrac{1}{2}+C\,\gamma_n^*+r_n+O(n^{-1/2})=\tfrac{1}{2}+C\,\gamma_n^*+O(n^{-1/2}),\end{aligned}$$

since $r_n=o(n^{-1/2})$. This proves (22). □

Theorem 5.2 follows from Lemma D.4.

## Appendix E: Missing proofs from Section 6

### E.1. Proof of Proposition 6.1

*Proof.* For $j\in\{0,1\}$, let

$$\widehat{S}(d_j)=\sum_i\frac{I(D_i=d_j)Y_i(D_i)}{\pi_i(d_j|n_{d_j})}$$

Then taking conditional expectation with respect to $n_{d_j}$, we get

$$E[\hat{S}(d_j)|n_{d_j}]=\sum_i\frac{E[I(D_i=d_j)|n_{d_j}]Y_i(d_j)}{\pi_i(d_j|n_{d_j})}=S(d_j).$$

Hence, $E[\hat{S}(d_j)]=S_i(d_j)$. This shows that $E[\hat{\theta}_{cht}]=\theta(d_1,d_0)$ □

### E.2. Proof of Theorem 6.1

*Proof.* Consider the following set of potential outcomes $Y_i^*(d_j) = \pi_i(d_j|n_{d_j})$. By the definition of the conditional H-T estimator,

$$\hat{\theta}_{CHT}(d_j) = \frac{1}{n}\sum_i \frac{I(D_i = d_j)}{\pi_i(d_j|n_{d_j})} Y_i^*(d_j) = \frac{1}{n}\sum_i I(D_i = d_j) = \frac{n_{d_j}}{n}$$

Next, we have

$$\theta^*(d_j) = \frac{1}{n}\sum_i \pi_i(d_j|n_{d_j}) = \frac{n_{d_j}}{n}$$

where the last equality follows from Lemma E.5. Hence at $Y_i^*(d_j)$, $\theta_{CHT}(d_j) = \theta^*(d_j)$, and $MSE(\hat{\theta}_{CHT}, \theta^*) = 0$ Also note that the weights $w_i(\mathbf{Z}, d_j) = \frac{1}{\pi_i(d_j|n_{d_j})} > 1$ as $0 < \pi_i(d_j|d_{d_j}) < 1$. Hence the conditions of Lemma 5.1 hold, which implies $\hat{\theta}_{CHT}$ is admissible. □

### E.3. Proof of Lemma 6.1

We begin with a few useful lemmas.

**Lemma E.1.** *Let $U$ and $V$ be two discrete random variables. Let $(U, V)$ denote the maximal coupling between them. Then $|P(U = t) - P(V = t)| \leq TV(U, V)$.*

*Proof.* Recall that for any two random variables $U$ and $V$, the maximal coupling always exists and is such that $P(U \neq V) = TV(U, V)$. Let $(U, V)$ be a maximal coupling of discrete random variables $U$ and $V$. In the following set of inequalities, we make repeated use of the fact that two events $A$ and $B$, we have, $P(A) - P(B^c) \leq P(A \cap B) \leq P(A)$.

$$\begin{aligned} P(U = t) &= P(U = t, U = V) + P(U = t, U \neq V) \\ &\leq P(V = t) + P(U \neq V) \end{aligned}$$

Similarly, we have

$$\begin{aligned} P(U = t) &= P(U = t, U = V) + P(U = t, U \neq V) \\ &\geq P(V = t, U = V) \\ &\geq P(V = t) - P(U \neq V) \end{aligned}$$

Combining these two inequalities, we get $P(U = t) - P(V = t) \leq P(U \neq V) = TV(U, V)$ and $P(U = t) - P(V = t) \geq P(U \neq V) = TV(U, V)$. Thus, $|P(U = t) - P(V = t)| \leq TV(U, V)$.

□

**Lemma E.2.** *Fix a treatment-exposure combination d, and let* $I_i = I(D_i = d)$ *and* $\pi_i = P(I_i = 1)$. *Denote* $B_i = Z_{N_i}(2)$. *Let* $s_n = \max_i \sum_{j \in B_i} \pi_j$ *and* $\Lambda_n = \sum_i \pi_i$. *Assume (a)* $\Lambda_n \leq C$, *(b)* $s_n \to 0$ *and (c) there exists* $C > 0$ *such that for every* $j \in B_i$, $P(I_i = 1, I_j = 1), \leq C\pi_i\pi_j$ *then* $TV(\sum_i I_i, \sum_i W_i) \to 0$.

*Proof.* To prove this result, we need the following classical approximation theorem from Arratia et al. (1989) that states that a sequence of dependent Bernoulli random variables can be approximated by a Poisson distribution:

**Lemma E.3** (Theorem 1, Arratia et al. (1989))**.** *Let* $\{X_i : i \in [n]\}$ *be dependent Bernoulli random variables with* $p_i = P(X_i = 1)$. *Let* $p_{ij} = E[X_i X_j]$. *For each i, let* $B_i \subset [n]$ *be a dependence neighborhood of i, i.e.* $X_i$ *is independent of* $X_j$ *for all* $j \notin B_i$. *Let* $b_1 = \sum_{i \in [n]} \sum_{j \in B_i} p_i p_j$, $b_2 = \sum_{i \in [n]} \sum_{j \neq i \in B_i} p_{ij}$ *and* $b_3 = \sum_{i \in [n]} E[E[X_i - p_i | \sigma(X_j : j \notin B_i)]]$. *Let Y denote a Poisson random variable with mean* $\lambda = \sum_i p_i$. *Then,* $TV(\sum_i X_i, Y) \leq 2(b_1 + b_2 + b_3)$

We apply the above lemma to the dependent indicators $I_i(D_i = d)$. Recall that $D_i$ is a function of $Z_i$ and $Z_{N_i}$. For any $i$ and $j$, $D_i$ is dependent on $D_j$ only if nodes $i$ and $j$ share a common neighbor, or if $i$ and $j$ are directly connected. Thus, using the notation in Lemma E.3, if we let $B_i = Z_{N_i}(2)$, $D_i$ is independent of $\{D_j, j \notin B_i\}$. Hence $b_3 = 0$. Moreover, we have the following:

$$b_1 = \sum_{i=1}^{n} \pi_i \sum_{j \in B_i} \pi_j \leq \big(\max_i \sum_{j \in B_i} \pi_j\big) \sum_{i=1}^{n} \pi_i = \Lambda_n\, s_n \to 0 \text{ by assumption (a) and (b).}$$

Similarly, from (c), for each dependent pair $(i, j)$, $\mathbb{P}(I_i = 1, I_j = 1) \leq C\, \pi_i \pi_j$, hence

$$b_2 = \sum_i \sum_{j \in B_i \setminus \{i\}} \mathbb{P}(I_i = 1, I_j = 1) \leq C \sum_i \sum_{j \in B_i} \pi_i \pi_j = C\, b_1 \to 0.$$

Hence $b_1 + b_2 \to 0$, implying $TV(\sum_i I_i, W) \to 0$. Note that

$$s_n = \max_i \sum_{j \in B_i} \pi_j = \max_i |B_i| \bar{\pi}_i \leq \pi_{\max} \max_i |B_i| \to 0,$$

if $\max_i |B_i| = o(n)$ and $\pi_{\max} = O(1/n)$. □

We are now ready to present a proof of Lemma 6.1. By definition, the conditional propensity score can be written as

$$P\left(I_i = 1 | \sum_i I_i = t\right) = \frac{P\left(I_i = 1, \sum_{j \neq i} I_j = t - 1\right)}{P\left(\sum_i I_i = t\right)} \tag{23}$$

We will show that both the numerator and denominator above can be approximated by the corresponding estimates from the vector $(W_1, \ldots, W_n)$, where each $W_i$ is an independent Poisson random variable with mean $\pi_i$. Let $U = \sum_i I_i$ and $V = \sum_i W_i$ and consider the maximal coupling between $U$ and $V$. Lemma E.1 gives us, $P(\sum_i I_i = t) \geq P(\sum_i W_i = t) - TV(\sum_i I_i, \sum_i W_i)$. Similarly,

$$
\begin{aligned}
& P\left(I_i = 1, \sum_{j\neq i} I_j = t-1\right) \\
& \leq P\left(I_i = 1, \sum_{j\neq i} I_j = t-1, I_i = W_i, \sum_{j\neq i} I_j = \sum_{j\neq i} W_j\right) \\
& + P\left(I_i = 1, \sum_{j\neq i} I_j = t-1, I_i \neq W_i, \sum_{j\neq i} I_j \neq \sum_{j\neq i} W_j\right) \\
& \leq P\left(W_i = 1, \sum_{j\neq i} W_j = t-1\right) + P\left(I_i \neq W_i, \sum_{j\neq i} I_j \neq \sum_{j\neq i} W_j\right) \\
& \leq P\left(W_i = 1, \sum_{j\neq i} W_j = t-1\right) + P\left(\sum_{j\neq i} I_j \neq \sum_{j\neq i} W_j\right) \\
& \leq P\left(W_i = 1, \sum_{j\neq i} W_j = t-1\right) + TV\left(\sum_{j\neq i} I_j, \sum_{j\neq i} W_j\right) d
\end{aligned}
$$

Thus, we get,

$$
\begin{aligned}
& \frac{P\left(I_i = 1, \sum_{j\neq i} I_j = t-1\right)}{P\left(\sum_i I_i = t\right)} \\
& \leq \frac{P\left(W_i = 1, \sum_{j\neq i} W_j = t-1\right) + TV(\sum_{j\neq i} I_j, \sum_{j\neq i} W_j)}{P\left(\sum_i W_i = t\right) - TV(\sum_i I_i, \sum_i W_i)} \\
& \to \frac{P\left(W_i = 1, \sum_{j\neq i} W_j = t-1\right)}{P\left(\sum_i W_i = t\right)} = P\left(W_i = 1 | \sum_i W_i = t\right)
\end{aligned}
$$

if $TV(\sum_i I_i, \sum_j W_j) \to 0$ and $TV(\sum_{j\neq i} I_j, \sum_{j\neq i} W_j) \to 0$ for each $i$. Apply Lemma E.2 directly to the first TV, and apply Lemma E.2 with leaving the node $i$ out to the second TV completes the proof.

### E.4. Proof of Corollary 1

**Lemma E.4.** *Let* $W_1, \ldots, W_n$ *be independent* $\mathrm{Pois}(\lambda_1), \ldots, \mathrm{Pois}(\lambda_n)$ *and write* $\lambda = \sum_{j=1}^{n} \lambda_j$ *and* $q_i = \lambda_i/\lambda$. *Let* $S_W = \sum_{j=1}^{n} W_j$ *and fix an integer* $t \geq 1$ *with* $P(S_W = t) > 0$. *Then* $|P(W_i = 1|S_w = t) - tq_i| \to 0$ *as long as* $\max_i q_i \to 0$.

*Proof.* For independent Poisson$(\lambda_j)$ variables the joint law conditional on the sum $S_W = t$ is Multinomial$(t; q_1, \ldots, q_n)$ with $q_j = \lambda_j/\lambda$. Hence, for small $q_i$ and fixed $t$,

$$P(W_i = 1 \mid S_W = t) = tq_i(1 - q_i)^{t-1} = tq_i(1 + O(tq_i)) = tq_i(1 + o(1))$$

which proves the claim.

□

*Proof.* As per Lemma 6.1, the vector $I_1, \ldots, I_n$ can be approximated by a vector $W_1, \ldots, W_n$, where each $W_i$ has Poisson distribution with mean $\pi_i$, From Lemma E.4, we have $P(W_i = 1| \sum W_i = t) \approx t\pi_i/(\sum_i \pi_i) \approx \pi_i(\sum_i I_i)/(\sum_i \pi_i)$. Substituting this in the definition of the conditional HT estimator, we obtain the ratio estimator.

□

### E.5. Proof of Theorem 6.2

*Proof.* We begin with two useful Lemmas. First, we show that the sum of conditional propensity score of treatment-exposure status $d_j$ is equal to the number of units observed under that treatment-exposure status.

**Lemma E.5.** $\sum_i \pi_i(d_j|n_{d_j}) = n_{d_j}$

*Proof.* We illustrate the proof for $d_1$. Let $U_i = 1$ when $D_i = d_1$ and 0 otherwise. Then $n_{d_1} = \sum_i U_i$, and the conditional propensity score for treatment-exposure status $d_1$ is given by $\pi_i(d_1|n_{d_1}) = \mathbb{P}(U_i = 1| \sum_i U_i)$. Using the identity

$$1 = \frac{\sum_i U_i}{\sum_i U_i},$$

and taking conditional expectation with respect to $\sum_i U_i$ on both sides, we get

$$1 = \sum_i E\left[\frac{U_i}{\sum_i U_i} | \sum_i U_i\right]$$
$$= \frac{1}{\sum_i U_i}\sum_i E\left[U_i | \sum_i U_i\right]$$
$$= \frac{1}{\sum_i U_i}\sum_i \mathbb{P}\left(U_i = 1 | \sum_i U_i\right) \quad (24)$$

This implies $\sum_i \mathbb{P}(U_i = 1| \sum_i U_i) = \sum_i U_i$, as desired. □

Next, we show that if the propensity scores are close to each other, the error of approximating the conditional H-T estimator with the difference-in-means estimator is small.

**Lemma E.6.** *Let $U_i$ be defined as in the previous lemma. Let $T = \sum_{i=1}^n U_i$, $p_i := P(U_i = 1 \mid T)$, $\bar{p} := \frac{T}{n}$. Let*

$$\beta_n := \frac{1}{n}\sum_{i=1}^n |p_i - \bar{p}|.$$

*Assume bounded potential outcomes, i.e. $Y_i(z_i, e_i) \leq M$ for all $i$ and that $p_i \geq p_{\min} > 0$ for all $i$. Define the estimators*

$$\hat{S}_1 = \frac{1}{n}\sum_{i=1}^n \frac{Y_i I(U_i = 1)}{p_i}, \qquad \hat{S}_{1,a} = \frac{\sum_{i=1}^n Y_i I(U_i = 1)}{T}.$$

*If $T/n \geq c > 0$ with high probability, then*

$$|\hat{S}_1 - \hat{S}_{1,a}| = O_p(\beta_n)$$

*Proof.* We begin by rewriting the difference:

$$\hat{S}_1 - \hat{S}_{1,a} = \frac{1}{n}\sum_{i=1}^n Y_i I(U_i = 1)\Big(\frac{1}{p_i} - \frac{n}{T}\Big) = \frac{1}{T}\sum_{i=1}^n Y_i I(U_i = 1)\frac{p_i - \bar{p}}{p_i}.$$

Taking absolute values and using the triangle inequality gives

$$|\hat{S}_1 - \hat{S}_{1,a}| \leq \frac{1}{T}\sum_{i=1}^n |Y_i| I(U_i = 1)\frac{|p_i - \bar{p}|}{p_i}.$$

Since $|Y_i| \leq M$ and $p_i \geq p_{\min}$,

$$|\hat{S}_1 - \hat{S}_{1,a}| \leq \frac{M}{p_{\min}} \cdot \frac{1}{T} \sum_{i=1}^{n} I(U_i = 1)\, |p_i - \bar{p}|.$$

Bounding the selected-unit sum by the full sum,

$$\frac{1}{T} \sum_{i=1}^{n} I(U_i = 1)\, |p_i - \bar{p}| \leq \frac{1}{T} \sum_{i=1}^{n} |p_i - \bar{p}| = \frac{n}{T} \beta_n,$$

and substituting this yields

$$|\hat{S}_1 - \hat{S}_{1,a}| \leq \frac{M}{p_{\min}} \cdot \frac{n}{T} \beta_n,$$

which proves the claim. □

Lemma E.5 implies that when the conditional propensity scores are equal to each other, we have $p_i = P(U_i | \sum_i U) = \sum_i U_i/n$. Substituting this in the definition of the conditional HT estimator shows that the conditional HT estimator reduces to the difference-in-means estimator. The approximation result follows by applying Lemma E.6, as follows. Let $\hat{S}_j = (1/n) \cdot \sum_i Y_i(d_j)(U_i/p_i)$ and $\hat{S}_{1,a} = \sum_i Y_i(d_j)U_i / \sum_i U_i$. Then, we have, $\Delta = \hat{\theta}_{CHT} - \hat{\theta}_{df} = \hat{S}_1 - \hat{S}_2 - (\hat{S}_{1,a} - \hat{S}_{2,a})$, so $|\Delta| \leq |\hat{S}_1 - \hat{S}_{1,a}| + |\hat{S}_2 - \hat{S}_{2,a}| = O(\beta_n)$ □

### *E.6. Proof of Theorem A.1*

*Proof.* For a CRD design, define,

$$\begin{aligned}
\alpha_i(1, e_i) &= \mathbb{E}\left[\frac{I(Z_i = 1, E_i = e_i)}{\sum_i Z_i}\right] \\
&= \frac{1}{n_t} \mathbb{P}\left(I(Z_i = 1, E_i = e_i)\right) \\
&= \frac{1}{n_t} \frac{n_t}{n} \frac{\binom{n_t - 1}{e_i}\binom{n_c}{d_i - e_i}}{\binom{n-1}{d_i}} \text{ if } n_t \geq e_i + 1 \text{ and } n_c \geq d_i - e_i, 0 \text{ otherwise} \\
&= \frac{1}{n} \frac{\binom{n_t - 1}{e_i}\binom{n_c}{d_i - e_i}}{\binom{n-1}{d_i}} \text{ if } n_t \geq e_i + 1 \text{ and } n_c \geq d_i - e_i, 0 \text{ otherwise}
\end{aligned}$$

For a Bernoulli trial, define,

$$
\begin{aligned}
\alpha_i(1, e_i) &= \mathbb{E}\left[\frac{I(Z_i = 1, E_i = e_i)}{\sum_i Z_i}\right] \\
&= \mathbb{E}_k\left[\mathbb{E}\left[\frac{I(Z_i = 1, E_i = e_i)}{\sum_i Z_i}\middle| \sum_i Z_i = k\right]\right] \\
&= \mathbb{E}\left[\frac{1}{\sum_i Z_i}\mathbb{P}\left(Z_i = 1, E_i = e_i \middle| \sum_i Z_i = k\right)\right] \\
&= \mathbb{E}_K\left[\frac{1}{K}\frac{K}{n}\frac{\binom{K-1}{e_i}\binom{n-K}{d_i-e_i}}{\binom{n-1}{d_i}}\right], \text{ where } \sum_i Z_i = K \\
&= \frac{1}{n}\mathbb{E}_K\left[\frac{\binom{K-1}{e_i}\binom{n-K}{d_i-e_i}}{\binom{n-1}{d_i}}\right]
\end{aligned}
$$

where $K$ is a restricted binomial random variable with support on $\{1, \ldots, N-1\}$ and $\mathbb{P}(K = k) = \frac{\binom{n}{k}p^k(1-p)^{n-k}}{1-(1-p)^n-p^n}$. A similar proof holds for the other cases. □

### E.7. Proof of Theorem A.2

*Proof.* We will illustrate the proof for one case, the rest are similar. Assume $n_c > d_i$, then

$$
\begin{aligned}
P(z_i = 1, e_i = 1) = P(z_i = 1)P(e_i = 1|z_i = 1) &= \frac{n_t}{n}\left[1 - P(e_i = 0|z_i = 1)\right] \\
&= \frac{n_t}{n}\left[1 - \frac{n_c}{n-1}\frac{n_c - 1}{n-2}\cdots\frac{n_c - d_i + 1}{n - d_i}\right] \\
&= \frac{n_t}{n}\left[1 - \frac{\binom{n_c}{d_i}}{\binom{n-1}{d_i}}\right]
\end{aligned}
$$

If $n_c < d_i$, then $P(e_i = 0|z_i = 1) = 0$. □

## References

Jürgen Angst and Guillaume Poly. A weak cramér condition and application to edgeworth expansions. *Electron. J. Probab*, 22(59):1–24, 2017.

Peter M Aronow. A general method for detecting interference between units in randomized experiments. *Sociological Methods & Research*, 41(1):3–16, 2012.

Peter M Aronow and Cyrus Samii. Estimating average causal effects under general interference, with application to a social network experiment. *Annals of Applied Statistics*, 11(4):1912–1947, 2017.

Richard Arratia, Larry Goldstein, and Louis Gordon. Two moments suffice for poisson approximations: the chen-stein method. *The Annals of Probability*, pages 9–25, 1989.

Susan Athey, Dean Eckles, and Guido W Imbens. Exact p-values for network interference. *Journal of the American Statistical Association*, pages 1–11, 2017.

Guillaume Basse and Avi Feller. Analyzing two-stage experiments in the presence of interference. *Journal of the American Statistical Association*, 113(521): 41–55, 2018.

Guillaume W Basse and Edoardo M Airoldi. Model-assisted design of experiments in the presence of network-correlated outcomes. *Biometrika*, 105(4):849–858, 2018.

Guillaume W Basse, Avi Feller, and Panos Toulis. Randomization tests of causal effects under interference. *Biometrika*, 106(2):487–494, 2019.

Jake Bowers, Mark M Fredrickson, and Costas Panagopoulos. Reasoning about interference between units: A general framework. *Political Analysis*, 21(1): 97–124, 2013.

Louis HY Chen and Qi-Man Shao. Normal approximation under local dependence. *Annals of probability*, 32(3A):1985–2028, 2004.

David Choi. New estimands for experiments with strong interference. *Journal of the American Statistical Association*, (just-accepted):1–22, 2023.

David Roxbee Cox. Planning of experiments. 1958.

Dean Eckles, Brian Karrer, and Johan Ugander. Design and analysis of experiments in networks: Reducing bias from interference. *Journal of Causal Inference*, 5 (1):20150021, 2016.

Mengsi Gao and Peng Ding. Causal inference in network experiments: regression-based analysis and design-based properties. *arXiv preprint arXiv:2309.07476*, 2023.

VP Godambe. A unified theory of sampling from finite populations. *Journal of the Royal Statistical Society. Series B (Methodological)*, pages 269–278, 1955.

VP Godambe and VM Joshi. Admissibility and bayes estimation in sampling finite populations. i. *The Annals of Mathematical Statistics*, 36(6):1707–1722, 1965.

J Hájek. Comment on a paper by d. basu. *Foundations of statistical inference*, 236, 1971.

Jeremy Heng, Pierre E Jacob, and Nianqiao Ju. A simple markov chain for independent bernoulli variables conditioned on their sum. *arXiv preprint arXiv:2012.03103*, 2020.

Daniel G Horvitz and Donovan J Thompson. A generalization of sampling without replacement from a finite universe. *Journal of the American statistical Association*, 47(260):663–685, 1952.

Yuchen Hu, Shuangning Li, and Stefan Wager. Average direct and indirect causal effects under interference. *Biometrika*, 109(4):1165–1172, 2022.

Michael G Hudgens and M Elizabeth Halloran. Toward causal inference with interference. *Journal of the American Statistical Association*, 103(482):832–842, 2008a.

Michael G Hudgens and M Elizabeth Halloran. Toward causal inference with interference. *Journal of the American Statistical Association*, (482):832–842, 2008b.

Ravi Jagadeesan, Natesh S Pillai, and Alexander Volfovsky. Designs for estimating the treatment effect in networks with interference. 2020.

Vishesh Karwa and Edoardo M Airoldi. A systematic investigation of classical causal inference strategies under mis-specification due to network interference. *arXiv preprint arXiv:1810.08259*, 2018.

Michael P Leung. Causal inference under approximate neighborhood interference. *Econometrica*, 90(1):267–293, 2022.

Shuangning Li and Stefan Wager. Random graph asymptotics for treatment effect estimation under network interference. *The Annals of Statistics*, 50(4):2334–2358, 2022.

Wenrui Li, Daniel L Sussman, and Eric D Kolaczyk. Causal inference under network interference with noise. *arXiv preprint arXiv:2105.04518*, 2021.

Charles F Manski. Identification of treatment response with social interactions. *The Econometrics Journal*, 16(1):S1–S23, 2013.

Kari Lock Morgan and Donald B Rubin. Rerandomization to improve covariate balance in experiments. *Annals of Statistics*, 40(2):1263–1282, 2012.

J Neyman. Sur les applications de la thar des probabilities aux experiences agaricales: Essay des principle. excerpts reprinted (1990) in english. *Statistical Science*, 5:463–472, 1923.

Jean Pouget-Abadie, Guillaume Saint-Jacques, Martin Saveski, Weitao Duan, Souvik Ghosh, Ya Xu, and Edoardo M Airoldi. Testing for arbitrary interference on experimentation platforms. *Biometrika*, 106(4):929–940, 2019.

James M Robins, Andrea Rotnitzky, and Lue Ping Zhao. Estimation of regression coefficients when some regressors are not always observed. *Journal of the*

*American statistical Association*, 89(427):846–866, 1994.

Paul R Rosenbaum. Reduced sensitivity to hidden bias at upper quantiles in observational studies with dilated treatment effects. *Biometrics*, 55(2):560–564, 1999.

Paul R Rosenbaum. Interference between units in randomized experiments. *Journal of the American Statistical Association*, 102(477), 2007.

Donald B Rubin. Estimating causal effects of treatments in randomized and nonrandomized studies. *Journal of Educational Psychology*, 66(5):688, 1974.

Donald B Rubin. Bayesian inference for causal effects: The role of randomization. *The Annals of statistics*, pages 34–58, 1978.

Fredrik Sävje, Peter Aronow, and Michael Hudgens. Average treatment effects in the presence of unknown interference. *Annals of statistics*, 49(2):673, 2021.

Michael E Sobel. What do randomized studies of housing mobility demonstrate? causal inference in the face of interference. *Journal of the American Statistical Association*, 101(476):1398–1407, 2006.

Daniel L Sussman and Edoardo M Airoldi. Elements of estimation theory for causal effects in the presence of network interference. *arXiv preprint arXiv:1702.03578*, 2017.

F Sävje. Causal inference with misspecified exposure mappings: separating definitions and assumptions. *Biometrika*, 03 2023. ISSN 1464-3510. . URL `https://doi.org/10.1093/biomet/asad019`. asad019.

Eric J Tchetgen Tchetgen and Tyler J VanderWeele. On causal inference in the presence of interference. *Statistical methods in medical research*, 21(1):55–75, 2012.

Panos Toulis and Edward Kao. Estimation of causal peer influence effects. In *International conference on machine learning*, pages 1489–1497. PMLR, 2013.

Johan Ugander and Hao Yin. Randomized graph cluster randomization. *Journal of Causal Inference*, 11(1):20220014, 2023.

Johan Ugander, Brian Karrer, Lars Backstrom, and Jon Kleinberg. Graph cluster randomization: Network exposure to multiple universes. In *Proceedings of the 19th ACM SIGKDD international conference on Knowledge discovery and data mining*, pages 329–337, 2013.

Davide Viviano. Experimental design under network interference. *arXiv preprint arXiv:2003.08421*, 2020.

Christina Lee Yu, Edoardo M Airoldi, Christian Borgs, and Jennifer T Chayes. Estimating the total treatment effect in randomized experiments with unknown network structure. *Proceedings of the National Academy of Sciences*, 119(44): e2208975119, 2022.